\documentclass[12pt]{article}
\usepackage{jheppub}
\pdfoutput=1
\usepackage{amsmath,bbm,array,amsfonts,graphicx,wrapfig,lscape,float,mathtools,multirow,longtable,amsthm}
\usepackage[dvipsnames]{xcolor}
\usepackage{stackrel}
\usepackage[all]{xy}
\usepackage{caption}
\usepackage{subcaption}
\usepackage[bottom]{footmisc}
\usepackage{mathrsfs}
\usepackage{tikz}
\usepackage{bm}
\usepackage{color}
\usepackage{enumerate}
\usetikzlibrary{decorations.pathreplacing}
\usepackage{diagbox}
\numberwithin{equation}{section}
\numberwithin{figure}{section}
\numberwithin{table}{section}
\usepackage{pgfplots}
\pgfplotsset{compat=1.14}
\usepackage{makecell}
\usepackage[utf8]{inputenc}
\usepackage[export]{adjustbox}[2011/08/13]
\usepackage{boldline}
\usepackage{comment}

\usepackage[autostyle=true]{csquotes}
\newtheorem{definition}{Definition}[section]
\newtheorem{theorem}{Theorem}[section]

\newtheorem{proposition}[theorem]{Proposition}
\newtheorem{conjecture}[theorem]{Conjecture} 

\newtheorem{example}{Example}

\usepackage{todonotes}

\usepackage[toc,page]{appendix}

\usepackage{tocloft}

\newcommand{\IP}{\mathbb{P}}

\newcommand{\IR}{\mathbb{R}}
\newcommand{\IC}{\mathbb{C}}

\newcommand{\IZ}{\mathbb{Z}}

	\title{Polytopes and Machine Learning}
	
	\author[a,b]{Jiakang Bao}
	\author[a,b,c,d]{Yang-Hui He}
	\author[a,b]{Edward Hirst}
	\author[e]{Johannes Hofscheier}
	\author[e]{Alexander Kasprzyk}
	\author[a]{Suvajit Majumder}

	\affiliation[a]{
		Department of Mathematics, City, University of London, EC1V 0HB, UK}
	\affiliation[b]{
	    London Institute for Mathematical Sciences, Royal Institution, London W1S 4BS, UK}
	\affiliation[c]{
		Merton College, University of Oxford, OX1 4JD, UK}
	\affiliation[d]{
		School of Physics, NanKai University, Tianjin, 300071, P.R. China}
	\affiliation[e]{
		School of Mathematical Sciences, University of Nottingham, Nottingham, NG7 2RD, UK}
	
	\emailAdd{jiakang.bao@city.ac.uk}
	\emailAdd{hey@maths.ox.ac.uk}
	\emailAdd{edward.hirst@city.ac.uk}
	\emailAdd{johannes.hofscheier@nottingham.ac.uk}
	\emailAdd{a.m.kasprzyk@nottingham.ac.uk}
	\emailAdd{suvajit.majumder@city.ac.uk}
	
	\preprint{
		\begin{flushright}
			LIMS-2021-011
		\end{flushright}
	}

	\abstract{We introduce machine learning methodology to the study of lattice polytopes. With supervised learning techniques, we predict standard properties such as volume, dual volume, reflexivity, etc, with accuracies up to 100\%. We focus on 2d polygons and 3d polytopes with Pl\"ucker coordinates as input, which out-perform the usual vertex representation.
	}

\begin{document}
	\maketitle

\section{Introduction}\label{intro}
Lattice polytopes -- convex bodies with integral vertices -- lie in the intersection of combinatorics, optimisation theory, number theory, geometry, and theoretical physics.
In combinatorics, problems include finding lattice points inside a lattice polytope (for example,~ \cite{ehrhart1967probleme,Betke1985LatticePI}).
In number theory, lattice polytopes constitute the cornerstone of the geometry of numbers (for example,~\cite{minkowski1910geometrie,poonen2000lattice}).
In geometry, lattice polytopes play a key role in the study of toric varieties, whose geometry can be described combinatorially in terms of cones and fans (for example,~\cite{cox2011toric,kasprzyk2011reflexive,batyrev2010generalization,batyrev2013lattice}).
In string theory, Batyrev--Borisov mirror symmetry~\cite{BatyrevBorisov+2011+39+66} -- describing how Calabi--Yau manifolds can be naturally constructed from a reflexive polytope -- initiated a large-scale collaboration between fundamental physics and computational algebraic geometry, culminating in the classification of reflexive polytopes in dimensions three and four~\cite{Kreuzer:1995cd,Kreuzer:1998vb,Kreuzer:2000xy}.
Details of these, and many other, applications of lattice polytopes can be found in~ \cite{rebhan2013strings,hibi2019algebraic}.

There has been considerable recent interest in applying techniques from data science and machine learning~(ML) to the study of pure mathematical data. Whilst initially this arose in large part from the investigation of the string theory/algebraic geometry landscape~\cite{He:2017aed,Krefl:2017yox,Ruehle:2017mzq,Carifio:2017bov,He:2017set,Bao:2020sqg}, it is natural to ask whether this paradigm can be applied to different disciplines within pure mathematics~\cite{He:2018jtw,He:2021oav,Bao:2021auj,He:2020eva,Bao:2020nbi,Bao:2021olg}.

Combinatorial geometry lends itself to~ML. The central objects -- lattice polytopes -- can be cast into a matrix of integers and fed to a neural network~(NN) classifier or regressor. We will discuss different representations of this input data in this paper. Meanwhile, much of the work within the geometry of polytopes involves the extraction of numerical invariants such as volume, integral point enumeration, identifying topological invariants of the toric variety, etc. Can~NN techniques ``learn'' and extrapolate from this underlying structure? We shall see that this is indeed possible.

We should point out that our point d'appui is \emph{not} ML but is how to use ML to study the mathematics.
Polytope ML is an interesting field in itself, where learning algorithms are developed to find the hyperplanes that define the convex body, q.v., \cite{gottlieb2018learning,polytopeMLtalk}.
What we shall hope for, instead, is whether standard ML algorithms can find patterns in lattice polytopes and related algebraic geometry that are difficult or unknown.

The paper is organized as follows. In \S\ref{preliminaries}, we give a brief review on polytopes. In \S\ref{datasets} the features of the data for polytopes used in machine learning are described. Then, in \S\ref{ml}, we present our main results on machine learning for polytopes in different dimensions. In particular, we predict the volumes\footnote{Note that unlike the nice paper \cite{Krefl:2017yox} where the authors tried to learn the associated Calabi-Yau volume, we are learning the polytope volume here.} and dual volumes from their Pl\"ucker coordinates. Other invariants, such as the Gorenstein index and codimension, are also examined.

One is referred to \url{https://github.com/edhirst/PolytopeML.git} for the data and $\texttt{Python}$ code scripts used in this paper.

\section{Polytope Preliminaries}\label{preliminaries}
We begin with a brief reminder of the relevant concepts, which also serves to lay out notation and conventions. See, for example, \cite{sommerville2020introduction} for more details.

Polytopes are geometric objects formed from $n$ points in $d$-dimensional Euclidean space. 
Specifically, 
\begin{definition}\label{def:P}
A convex polytope $P$  has two equivalent definitions:
\begin{description}
    \item[Vertex Representation:] 
    the convex hull of set of $n$ points (called vertices)
      $p_i \in \IR^d$
\begin{equation}
    {\rm Conv}(\{p_i\}) = \left\{
      \sum\limits_{i=1}^n \alpha_i p_i \,\Big|\,
      \alpha_i \ge 0, \ \sum\limits_{i=1}^n \alpha_i = 1
      \right\};
\end{equation}
    \item[Half-hyperplane Representation:]
      the intersection of linear inequalities (hyperplanes called facets)
      \begin{equation}
          H \cdot \underline{x} \geq \underline{b} \ ,
      \end{equation}
      where $\underline{b}$ and $\underline{x}$ are real $n,d$-vectors respectively, and $H$ is some $n \times d$ matrix. 
\end{description}
The vertices are considered as 0-faces, the extremal edges are called 1-faces, and so on, until codimension 1, which are the $(n-1)$-faces that are also called facets.
The simplest polytope in $\IR^d$, consisting of $d+1$ vertices, is called a \textbf{simplex}. 
\end{definition}
In this paper, we will exclusively deal with \textbf{lattice polytopes}
where the vertices are integral points (lattice vectors), i.e., all $p_i \in \IZ^d$.
Henceforth, by polytope, we will mean a (convex) lattice polytope, denoted by $P$.

Given polytope $P$, a key concept is \emph{duality}:
\begin{definition}
The (polar) dual of $P$ is the polyhedron
\begin{equation}
    P^\circ \vcentcolon = \{\textbf{v} \in \mathbb{R}^d \; | \; \textbf{u}\cdot\textbf{v} \geq -1 \; \forall \; \textbf{u} \in P \}.
\end{equation}
\end{definition}
When $\textbf{0}\in\textup{int}(P)$ we have that $P^\circ$ is a polytope, albeit with rational-valued vertices.
Another important action on $P$ is:
\begin{definition}
The $r$-dilation of $P$ for $r \in \IZ_{\geq0}$ is the polytope
\begin{equation}
    rP := \{r \textbf{x} \; | \; \textbf{x} \in P \}.
\end{equation}
\end{definition}
Finally, we need a notion of primitivity:
\begin{definition}
A vector with integer components is called \textbf{primitive} if the GCD of the components is equal to $\pm 1$.
\end{definition}

In this paper, we are concerned with the following types of polytopes.
\begin{definition}
For a lattice polytope $P$, we have the following:
\begin{description}
\item[Fano: ] A convex polytope $P$ is called \emph{Fano} if the origin is strictly in its interior $\textup{int}(P)$ and all the vertices are primitive lattice points
(q.v.\cite{Kasprzyk_2010,Kasprzyk2012FANOP,balletti2016maximum,batyrev2019fine}).
These polytopes are called \emph{canonical Fano} if further $\textup{int}(P)\cap\mathbb{Z}^d=\{\mathbf{0}\}$.

\item[Reflexive: ] A Fano polytope $P$ is called \emph{reflexive} if its polar dual $P^\circ$ is also Fano. Indeed, for a generic lattice polytope the dual will have vertices in $\mathbb{Q}^d$, for reflexives, the vertices of the dual are also in $\IZ^d$.

\end{description}
\end{definition}
The terminology clearly comes from the study of toric geometry and these terms also apply to their corresponding toric varieties.
Importantly, we have that \cite{borisov1992singular,borisov2000convex}
\begin{theorem}
In any dimension $d$, there is a finite number of canonical Fano polytopes, as well as reflexive polytopes, up to $\textup{GL}(d, \IZ)$, which acts on the vertices.
\end{theorem}

\begin{example}
It is illustrative to give a concrete running example.
Consider the lattice polygon in Figure \ref{polygonexample}(a).
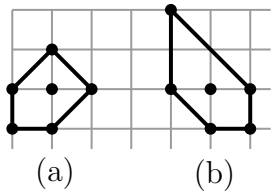
\begin{figure}[h]
    \centering
\tikzset{every picture/.style={line width=0.75pt}} 
\begin{tikzpicture}[x=0.75pt,y=0.75pt,yscale=-1,xscale=1]
\draw  [draw opacity=0] (92,42) -- (222.3,42) -- (222.3,112.2) -- (92,112.2) -- cycle ; \draw  [color={rgb, 255:red, 155; green, 155; blue, 155 }  ,draw opacity=1 ] (92,42) -- (92,112.2)(112,42) -- (112,112.2)(132,42) -- (132,112.2)(152,42) -- (152,112.2)(172,42) -- (172,112.2)(192,42) -- (192,112.2)(212,42) -- (212,112.2) ; \draw  [color={rgb, 255:red, 155; green, 155; blue, 155 }  ,draw opacity=1 ] (92,42) -- (222.3,42)(92,62) -- (222.3,62)(92,82) -- (222.3,82)(92,102) -- (222.3,102) ; \draw  [color={rgb, 255:red, 155; green, 155; blue, 155 }  ,draw opacity=1 ]  ;
\draw [line width=1.5]    (112,62) -- (132,82) ;
\draw  [fill={rgb, 255:red, 0; green, 0; blue, 0 }  ,fill opacity=1 ] (114.5,82) .. controls (114.5,80.62) and (113.38,79.5) .. (112,79.5) .. controls (110.62,79.5) and (109.5,80.62) .. (109.5,82) .. controls (109.5,83.38) and (110.62,84.5) .. (112,84.5) .. controls (113.38,84.5) and (114.5,83.38) .. (114.5,82) -- cycle ;
\draw [line width=1.5]    (112,62) -- (92,82) ;
\draw [line width=1.5]    (132,82) -- (112,102) ;
\draw [line width=1.5]    (92,82) -- (92,102) ;
\draw [line width=1.5]    (92,102) -- (112,102) ;
\draw  [fill={rgb, 255:red, 0; green, 0; blue, 0 }  ,fill opacity=1 ] (114.5,62) .. controls (114.5,60.62) and (113.38,59.5) .. (112,59.5) .. controls (110.62,59.5) and (109.5,60.62) .. (109.5,62) .. controls (109.5,63.38) and (110.62,64.5) .. (112,64.5) .. controls (113.38,64.5) and (114.5,63.38) .. (114.5,62) -- cycle ;
\draw  [fill={rgb, 255:red, 0; green, 0; blue, 0 }  ,fill opacity=1 ] (134.5,82) .. controls (134.5,80.62) and (133.38,79.5) .. (132,79.5) .. controls (130.62,79.5) and (129.5,80.62) .. (129.5,82) .. controls (129.5,83.38) and (130.62,84.5) .. (132,84.5) .. controls (133.38,84.5) and (134.5,83.38) .. (134.5,82) -- cycle ;
\draw  [fill={rgb, 255:red, 0; green, 0; blue, 0 }  ,fill opacity=1 ] (114.5,102) .. controls (114.5,100.62) and (113.38,99.5) .. (112,99.5) .. controls (110.62,99.5) and (109.5,100.62) .. (109.5,102) .. controls (109.5,103.38) and (110.62,104.5) .. (112,104.5) .. controls (113.38,104.5) and (114.5,103.38) .. (114.5,102) -- cycle ;
\draw  [fill={rgb, 255:red, 0; green, 0; blue, 0 }  ,fill opacity=1 ] (94.5,102) .. controls (94.5,100.62) and (93.38,99.5) .. (92,99.5) .. controls (90.62,99.5) and (89.5,100.62) .. (89.5,102) .. controls (89.5,103.38) and (90.62,104.5) .. (92,104.5) .. controls (93.38,104.5) and (94.5,103.38) .. (94.5,102) -- cycle ;
\draw  [fill={rgb, 255:red, 0; green, 0; blue, 0 }  ,fill opacity=1 ] (94.5,82) .. controls (94.5,80.62) and (93.38,79.5) .. (92,79.5) .. controls (90.62,79.5) and (89.5,80.62) .. (89.5,82) .. controls (89.5,83.38) and (90.62,84.5) .. (92,84.5) .. controls (93.38,84.5) and (94.5,83.38) .. (94.5,82) -- cycle ;
\draw [line width=1.5]    (192,102) -- (212,102) ;
\draw [line width=1.5]    (172,82) -- (192,102) ;
\draw [line width=1.5]    (172,42) -- (172,82) ;
\draw [line width=1.5]    (172,42) -- (212,82) ;
\draw [line width=1.5]    (212,82) -- (212,102) ;
\draw  [fill={rgb, 255:red, 0; green, 0; blue, 0 }  ,fill opacity=1 ] (194.5,82) .. controls (194.5,80.62) and (193.38,79.5) .. (192,79.5) .. controls (190.62,79.5) and (189.5,80.62) .. (189.5,82) .. controls (189.5,83.38) and (190.62,84.5) .. (192,84.5) .. controls (193.38,84.5) and (194.5,83.38) .. (194.5,82) -- cycle ;
\draw  [fill={rgb, 255:red, 0; green, 0; blue, 0 }  ,fill opacity=1 ] (174.5,82) .. controls (174.5,80.62) and (173.38,79.5) .. (172,79.5) .. controls (170.62,79.5) and (169.5,80.62) .. (169.5,82) .. controls (169.5,83.38) and (170.62,84.5) .. (172,84.5) .. controls (173.38,84.5) and (174.5,83.38) .. (174.5,82) -- cycle ;
\draw  [fill={rgb, 255:red, 0; green, 0; blue, 0 }  ,fill opacity=1 ] (174.5,42) .. controls (174.5,40.62) and (173.38,39.5) .. (172,39.5) .. controls (170.62,39.5) and (169.5,40.62) .. (169.5,42) .. controls (169.5,43.38) and (170.62,44.5) .. (172,44.5) .. controls (173.38,44.5) and (174.5,43.38) .. (174.5,42) -- cycle ;
\draw  [fill={rgb, 255:red, 0; green, 0; blue, 0 }  ,fill opacity=1 ] (194.5,102) .. controls (194.5,100.62) and (193.38,99.5) .. (192,99.5) .. controls (190.62,99.5) and (189.5,100.62) .. (189.5,102) .. controls (189.5,103.38) and (190.62,104.5) .. (192,104.5) .. controls (193.38,104.5) and (194.5,103.38) .. (194.5,102) -- cycle ;
\draw  [fill={rgb, 255:red, 0; green, 0; blue, 0 }  ,fill opacity=1 ] (214.5,102) .. controls (214.5,100.62) and (213.38,99.5) .. (212,99.5) .. controls (210.62,99.5) and (209.5,100.62) .. (209.5,102) .. controls (209.5,103.38) and (210.62,104.5) .. (212,104.5) .. controls (213.38,104.5) and (214.5,103.38) .. (214.5,102) -- cycle ;
\draw  [fill={rgb, 255:red, 0; green, 0; blue, 0 }  ,fill opacity=1 ] (214.5,82) .. controls (214.5,80.62) and (213.38,79.5) .. (212,79.5) .. controls (210.62,79.5) and (209.5,80.62) .. (209.5,82) .. controls (209.5,83.38) and (210.62,84.5) .. (212,84.5) .. controls (213.38,84.5) and (214.5,83.38) .. (214.5,82) -- cycle ;

\draw (102,115) node [anchor=north west][inner sep=0.75pt]   [align=left] {(a)};
\draw (182,116) node [anchor=north west][inner sep=0.75pt]   [align=left] {(b)};
\end{tikzpicture}
    \caption{The lattice polygons which are dual of each other.}\label{polygonexample}
\end{figure}
The vertices are $v_1=(1,0),v_2=(0,-1),v_3=(-1,-1),v_4=(-1,0),v_5=(0,1)$. 
One can check that its dual is the polygon shown in Figure \ref{polygonexample}(b) whose vertices are $v'_1=(1,0),v'_2=(1,-1),v'_3=(0,-1),v'_4=(-1,0),v'_5=(-1,2)$. 
Indeed, they are both lattice polygons, having vertices in $\mathbb{Z}^2$. 
Because $P$ and $P^\circ$ here are both lattice, this is a pair of reflexive polytopes in dimension 2.
Furthermore, they are clearly both canonical Fano.
\end{example}

We will also investigate the efficiency of ML in predicting the codimension of the embedding of the toric variety into the weighted projective space induced by the given polytope.
For instance, let $P \subset \IR^d$ be a Fano polytope and let $P^\circ$ its polar dual.
Lattice points in the cone over the polytope $P^\circ \times \{ 1 \} \subset \IR^{d+1}$ yield a semigroup which  accepts a unique Hilbert basis.
Let $n$ be the cardinality of this Hilbert basis.
The \emph{codimension} of $P$ is $n-\dim(P)-1$.

An important property of $P$ is its volume (which is clearly a $\textup{GL}(d, \IZ)$-invariant quantity).
For a simplex, there is a standard formula (q.v.~\cite{sommerville2020introduction})
\begin{theorem}[Cayley-Menger]
The volume of simplex in $\IR^d$ is 
$V_d = \frac{1}{d!} \prod\limits_{i=1}^d h_i$,
where $h_1$ is the distance between the first two vertices (the final answer is independent of the ordering of the vertices), $h_2$ is the height of the third vertex above the line containing the first two, $h_3$ is the height of the fourth vertex above the plane containing the first three vertices, etc.
\end{theorem}
For general polytope, we first have that:
\begin{definition}
A \textbf{triangulation} of a polytope $P$ is a partition of $P$ into \emph{simplices} such that
(i) the union of all them equals $P$ and (ii) the intersection of any pair of them is a (possibly empty) common face.
\end{definition}
In other words, a triangulation is a covering of $P$ by geometrical simplices.
Thereby, we have
\begin{definition}
The \textbf{volume} of a polytope $P$ is the sum of volumes of the simplices which constitutes any triangulation.
\end{definition}
The \emph{normalized volume} of a $d$-dimensional polytope is given by $d!$ times the volume.
We will also be interested in the volume of the dual polytope $P^\circ$, to which we will refer as \textbf{dual volume}.

Finally we introduce the Gorenstein index of $P$:
\begin{definition}
Let $P$ be a Fano polytope. The dual polytope $P^\circ$ is, in general, a rational polytope. Let $r\in\IZ_{>0}$ be the smallest positive integer such that $rP^\circ$ is a lattice polytope (and hence Fano). We call $r$ the \emph{Gorenstein index} of $P$.
\end{definition}
Clearly a Fano polytope is reflexive if and only if the Gorenstein index is $1$.

\begin{example}
For our running example in Figure \ref{polygonexample}(a), the polygon has (normalized) volume 5. Its dual volume, which is the volume of the polygon in Figure \ref{polygonexample}(b), is equal to 7. Since in general the dual polygon can have rational vertices (instead of only integer ones), the dual volume can be a rational number. The Gorenstein index is 1 in both cases, since both polygons are reflexive. The codimension of $P$ is $8-2-1=5$, and the codimension of $P^\circ$ is $6-2-1=3$.
\end{example}

\subsection{Pl\"ucker Coordinates}
In Definition \ref{def:P}, we gave two standard representations of $P$.
However, there is another important representation which we will use in this paper: the Pl\"ucker coordinate representation.

Many problems in mathematics and physics rely on the structure of polytopes, without regard for the specific embedding in the lattice they are drawn on.
It is with this motivation that Pl\"ucker coordinates are a natural consideration for polytope representation in our subsequent work.

\begin{definition}
Let $P$ be an $d$-dimensional polytope with $n$ vertices.
This structure can be reformatted into an $d \times n$ matrix $V$ where each column corresponds to a vertex $v_i$.
Consider the integer kernel of $V$ (aka grading) and take the maximal minors of $\ker(V)$;
this gives a list of integers known as the Pl\"ucker coordinate, which we regard as a point in projective space.
\end{definition}

\begin{example}
For the pentagon in Figure \ref{polygonexample}(a), the vertices matrix and its kernel are
\begin{equation}
V = 
\begin{pmatrix}
    v_1 & v_2 & v_3 & v_4 & v_5 \\ \hline
    1 & 0 &  -1 & -1 & 0 \\
    0 & -1 & -1 & 0 & 1 \\
\end{pmatrix}
\Rightarrow
\ker(V) =
    \begin{pmatrix}
        v_1 & v_2 & v_3 & v_4 & v_5 \\ \hline
        1 & 0 & 1 & 0 & 1 \\
        1 & 0 & 0 & 1 & 0 \\
        0 & 1 & 0 & 0 & 1 \\
    \end{pmatrix},\label{gradingexample}
\end{equation}
which in fact reveals the linear relations among the vertices: $v_1+v_3+v_5=0$, $v_1+v_4=0$ and $v_2+v_5=0$. By taking all of the possible $3 \times 3$ minors, of which there are ${5 \choose 3} = 10$, we get the Pl\"ucker coordinate as $(1:-1:1:0:-1:1:1:0:-1:1)\in\IC\IP^9$. Likewise, the polygon in Figure \ref{polygonexample}(b) has Pl\"ucker coordinate $(2:-1:1:-1:-1:1:2:0:-1:1)\in\IC\IP^9$.
\end{example}

\paragraph{Remark:} Note that due to the GL$(d,\mathbb{Z})$ invariance, it is possible that two inequivalent lattice polytopes can have the same Pl\"ucker coordinates. 
For example, the polygon with vertices $\{(-1,-1), (1,0), (0,1)\}$ and that with
$\{(-1,-1), (2,-1), (-1,2)\}$ both have Pl\"ucker coordinates $(1:1:1)\in\IC\IP^2$.
We can distinguish such polytopes by using an additional piece of combinatorial data called the quotient gradings. 
However, we will not take this into consideration in this paper. Instead, \emph{we require from here onwardwards that the vertices of $P$ generate the lattice $\IZ^d$}. With this additional assumption, the Pl\"ucker coordinates uniquely determine $P$ (up to the action of GL$(d,\mathbb{Z})$).

There is some trivial redundancy in the Pl\"ucker coordinates, based on the initial ordering of the polytope's vertices (columns of the vertex matrix). 
This leads to at most $n!$ equivalent ways of writing the Pl\"ucker coordinates for a polytope with $n$ vertices, coming from the permutation group $\mathfrak{S}_n$ shuffling the vertices. Such ordering is irrelevant for the Pl\"ucker coordinate object, but for tensor (vector) representation for ML input it allows some redundancy for data augmentation, as we will see.

It is not hard to check that the Pl\"ucker coordinates are GL$(d,\mathbb{Z})$ invariant. As aforementioned, Pl\"ucker coordinates reveal the linear relations among the vertices and are hence endowed with geometric meanings. These properties make the Pl\"ucker coordinates favourable as input vectors for our machine learning problem.

We summarize the key properties of Pl\"ucker coordinates as follows.
\begin{proposition}\label{redundancy}
Given a $d$-dimensional polytope $P$ with $n$ vertices, the number of Pl\"ucker coordinates is $N={n \choose n-d}$ as the $n$ vertices of a $d$-dimensional polytope satisfy $(n-d)$ linear relations, with (at most) $n!$ choices of Pl\"ucker coordinates in $\mathbb{CP}^{N-1}$. The Pl\"ucker coordinates are \textup{GL}$(d,\mathbb{Z})$ invariant.
\end{proposition}

\section{Polytope Datasets}\label{datasets}
Having introduced the preliminaries of polytopes which will be used in the paper, we now move on to describe the data we considered.

\subsection{Polygons in dimension 2}
Here we consider $P$ a Fano polygon whose vertices generate the lattice $\IZ^2$.
The dataset of such polygons used is available from the link in \S\ref{intro}. Each polygon has Gorenstein index at most $30$.
The dataset consists of polygons with vertices $n \in \{3,4,5,6\}$ with frequencies  $\{277,7041,16637,3003\}$, respectively, giving a total of $26958$ polygons.

For each polygon we compute the series of invariants:
\begin{equation}\label{Plabels}
(\mbox{volume, 
dual volume, 
Gorenstein index, 
codimension})
\end{equation}
These properties take values within the ranges: $[3,514], (0.21,15.37), [1,30], [2,42]$ respectively, where all properties are integer-valued except dual volume which are rational and hence computationally represented with floating-point reals.

\begin{figure}[h]
	\centering
	\begin{subfigure}{6cm}
		\centering
		\includegraphics[width=6cm]{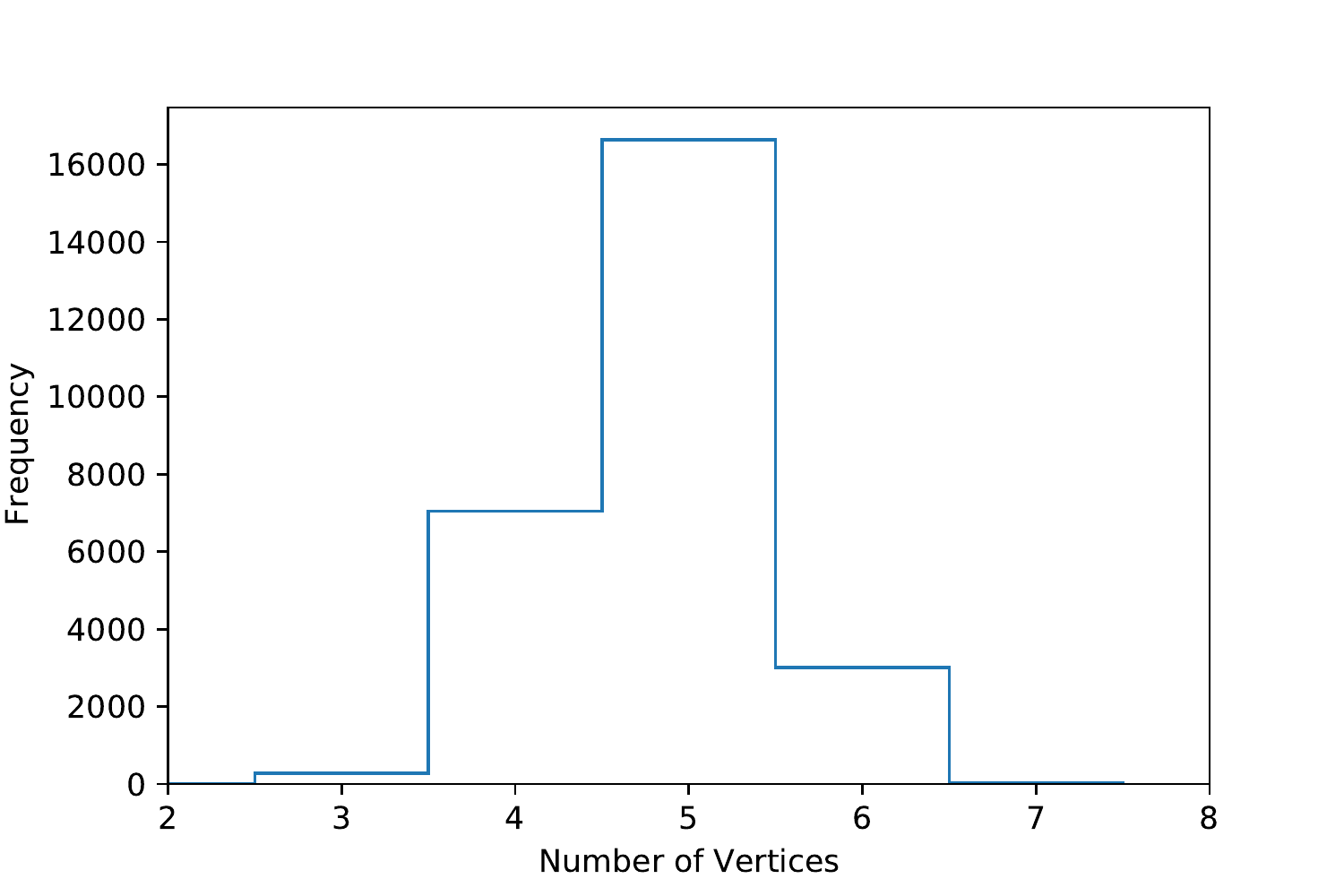}
		\caption{}\label{polygon_NoV}
	\end{subfigure} \\
    \begin{subfigure}{6cm}
    	\centering
    	\includegraphics[width=6cm]{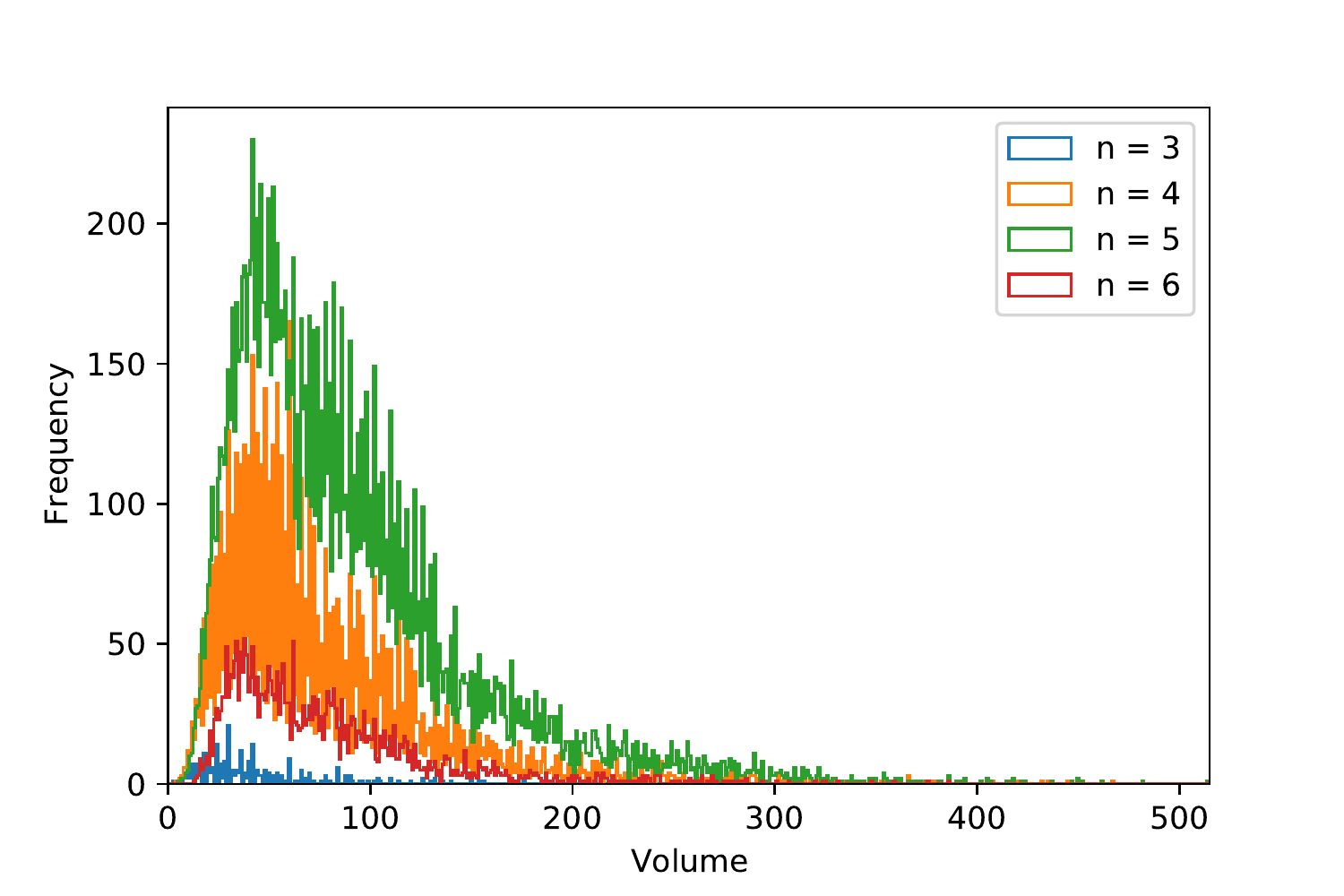}
    	\caption{}\label{polygon_Vol}
    \end{subfigure}
    \begin{subfigure}{6cm}
    	\centering
    	\includegraphics[width=6cm]{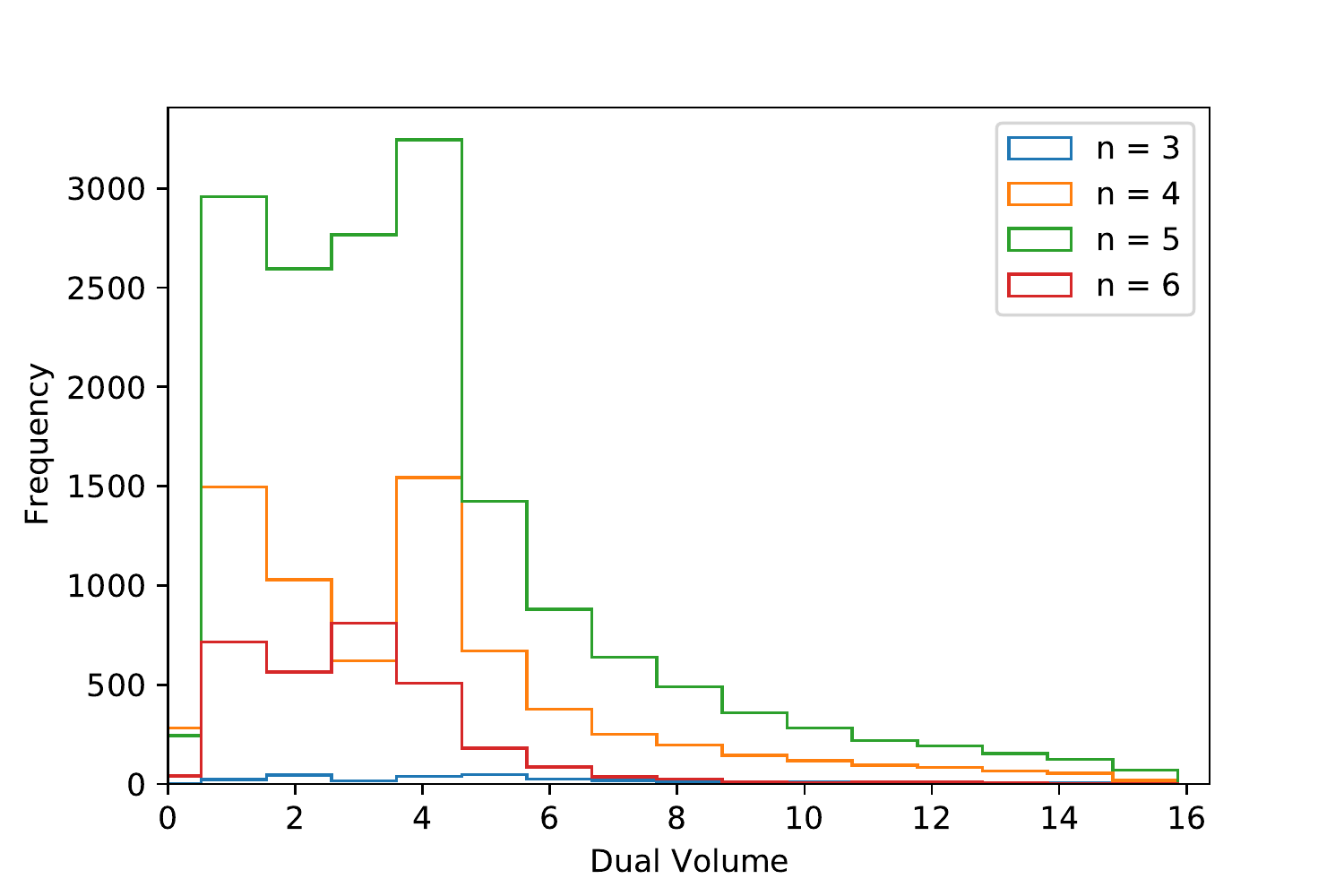}
    	\caption{}\label{polygon_DualVol}
    \end{subfigure} \\
    \begin{subfigure}{6cm}
		\centering
		\includegraphics[width=6cm]{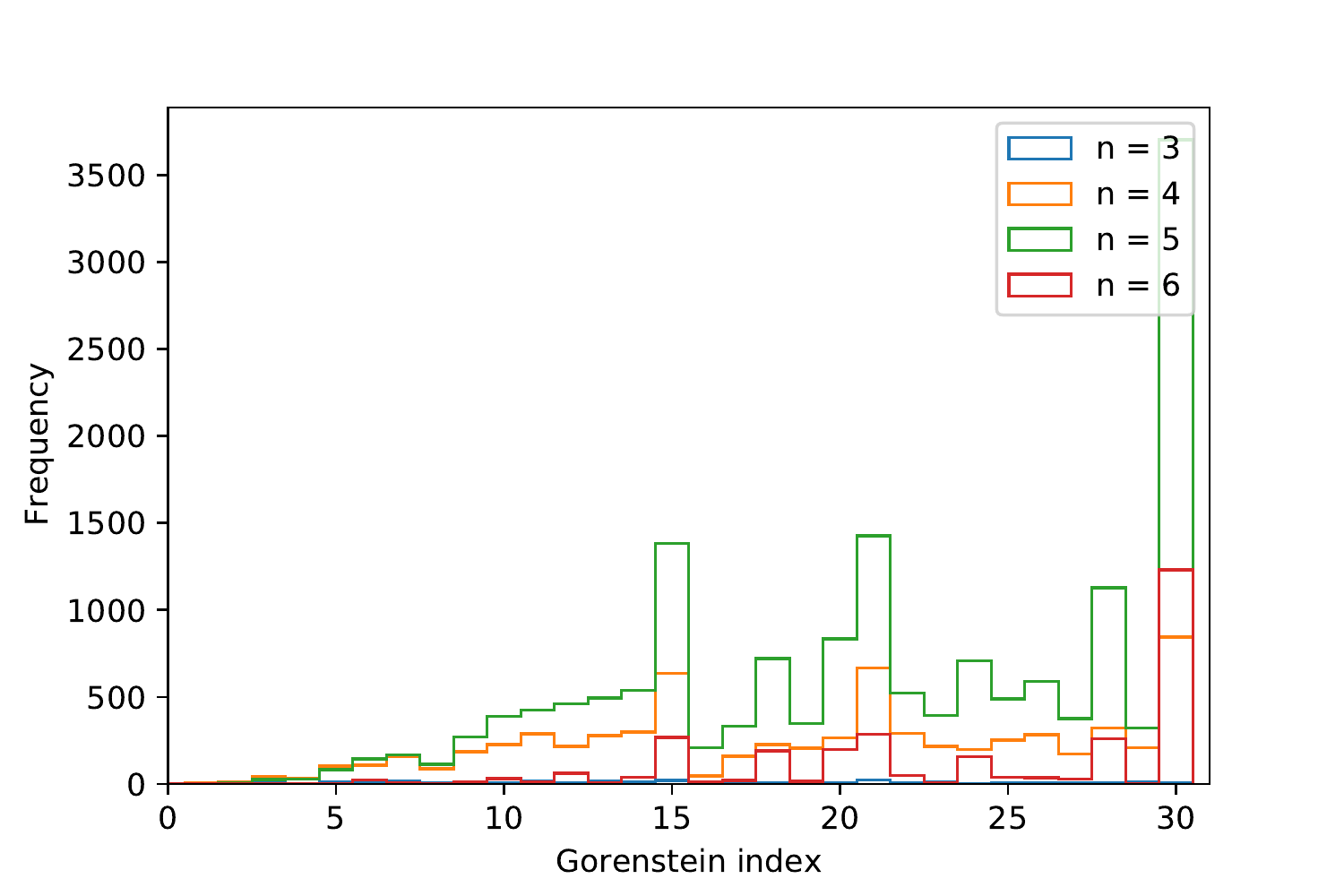}
		\caption{}\label{polygon_GI}
	\end{subfigure}
    \begin{subfigure}{6cm}
    	\centering
    	\includegraphics[width=6cm]{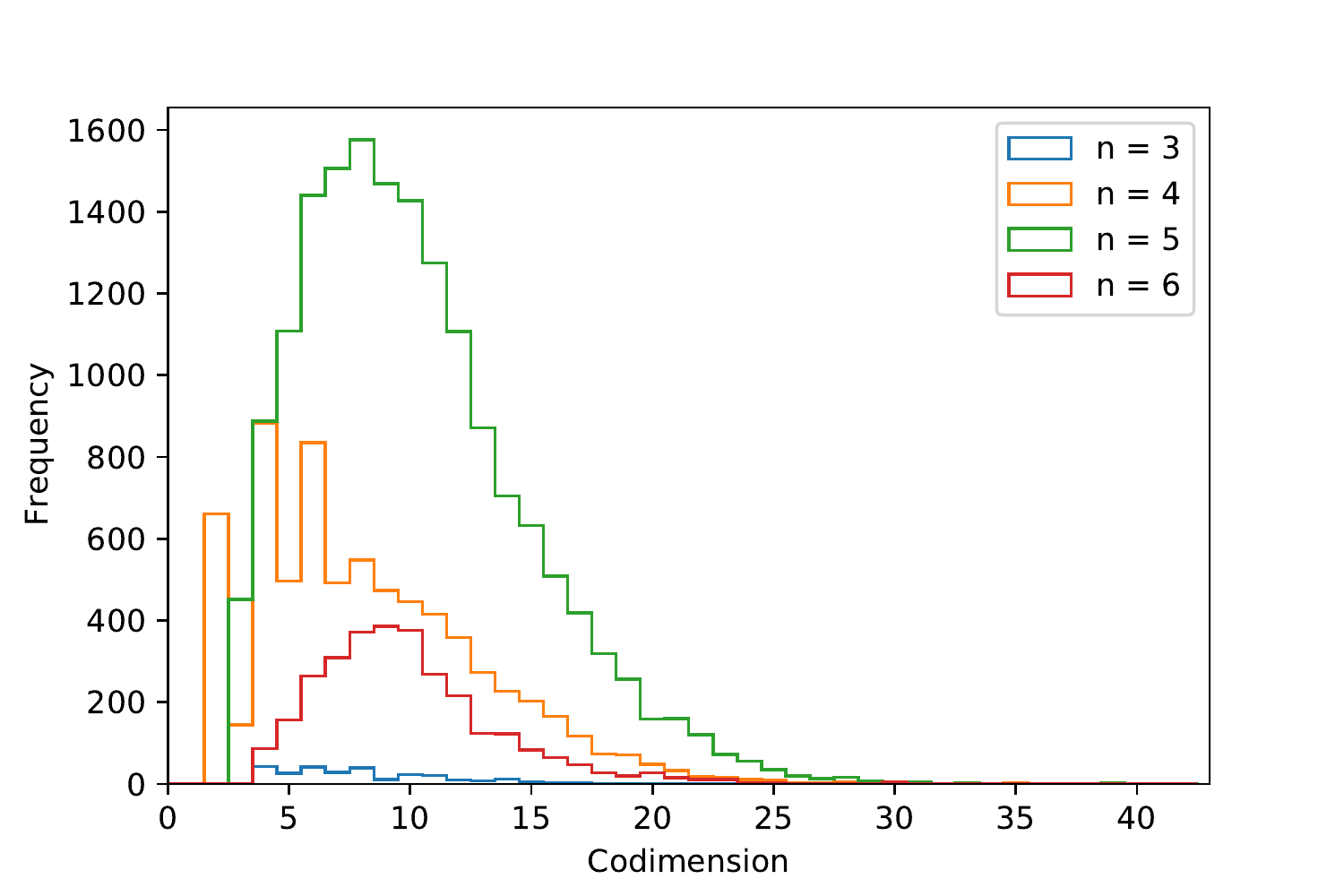}
    	\caption{}\label{polygon_Codim}
    \end{subfigure}
    \caption{Distributions of the polygon dataset's properties, split by the number of polygon's vertices, $n$. (a) shows the frequencies of polygons with each number of vertices in the dataset. (b-e) show the distributions over: volume, dual volume, Gorenstein index, and codimension respectively.}\label{polygon_dists}
\end{figure}

Plots of the frequencies of each of these properties across the dataset for each number of vertices are shown in Figure \ref{polygon_dists}.
Furthermore, Figure \ref{polygon_NoV} shows the distribution of the number of vertices across the dataset.

Now, we need to represent the polygons as appropriate tensors for the ML input. 
The naive initial choice of representation was a flattened list of the vertices, however as will be shown in the first investigation of \S\ref{polygon_ml} this is outperformed by the Pl\"ucker coordinate representation. 

Since for $n$ vertices the Pl\"ucker coordinate has length ${n \choose n-2}$, the input vectors vary in length.
Furthermore, there is redundancy in the ordering of the Pl\"ucker coordinate (or, more precisely, the ordering of the vertices when computing the Pl\"ucker coordinate), as discussed in Proposition \ref{redundancy}.
The dataset was thus augmented with (at most) 3 randomly chosen distinct Pl\"ucker representations for each polygon. This serves two roles: it prevents the machine from inadvertently learning in a way that depends on our choice of vertex ordering; and it provides a larger dataset on which to train.
This made the final dataset a list of 80874 Pl\"ucker representations with their respective polygon invariants.
For reference, the Pl\"ucker coordinates over the polygon dataset take value in the range [-450,450], whilst the coordinates of the vertices take value in the range [-290,180] (although this, of course, depends on the choice of basis).

\subsection{Polytopes in dimension 3}
Next, we move onto dimension 3.
Here our focus is the 674688 canonical Fano polytopes, as classified in \cite{Kas08}\footnote{The data can be found at \url{http://www.grdb.co.uk/forms/toricf3c} in the GRDB database.}. The full data contains many invariants for each polytope, including the volume, the dual volume, whether it is reflexive, the Gorenstein index, and the codimension.
As before, we enhance our dataset using the $n!$ choices of Pl\"ucker coordinate for a polytope with $n$ vertices. We randomly choose (at most) 10 distinct Pl\"ucker coordinates for each polytope. This augmentation of the data removes bias due to our choice of ordering of the vertices (and hence Pl\"ucker coordinate). The enhanced dataset  has 6744246 data points. 

Since this is a huge dataset, and since we would only need much fewer data points to train the model, we will consider 800000 random samples here. As the choice is completely random and this is still a large amount of data, we should lose no generality here.

\begin{figure}[h]
	\centering
	\begin{subfigure}{6cm}
	\includegraphics[width=6cm]{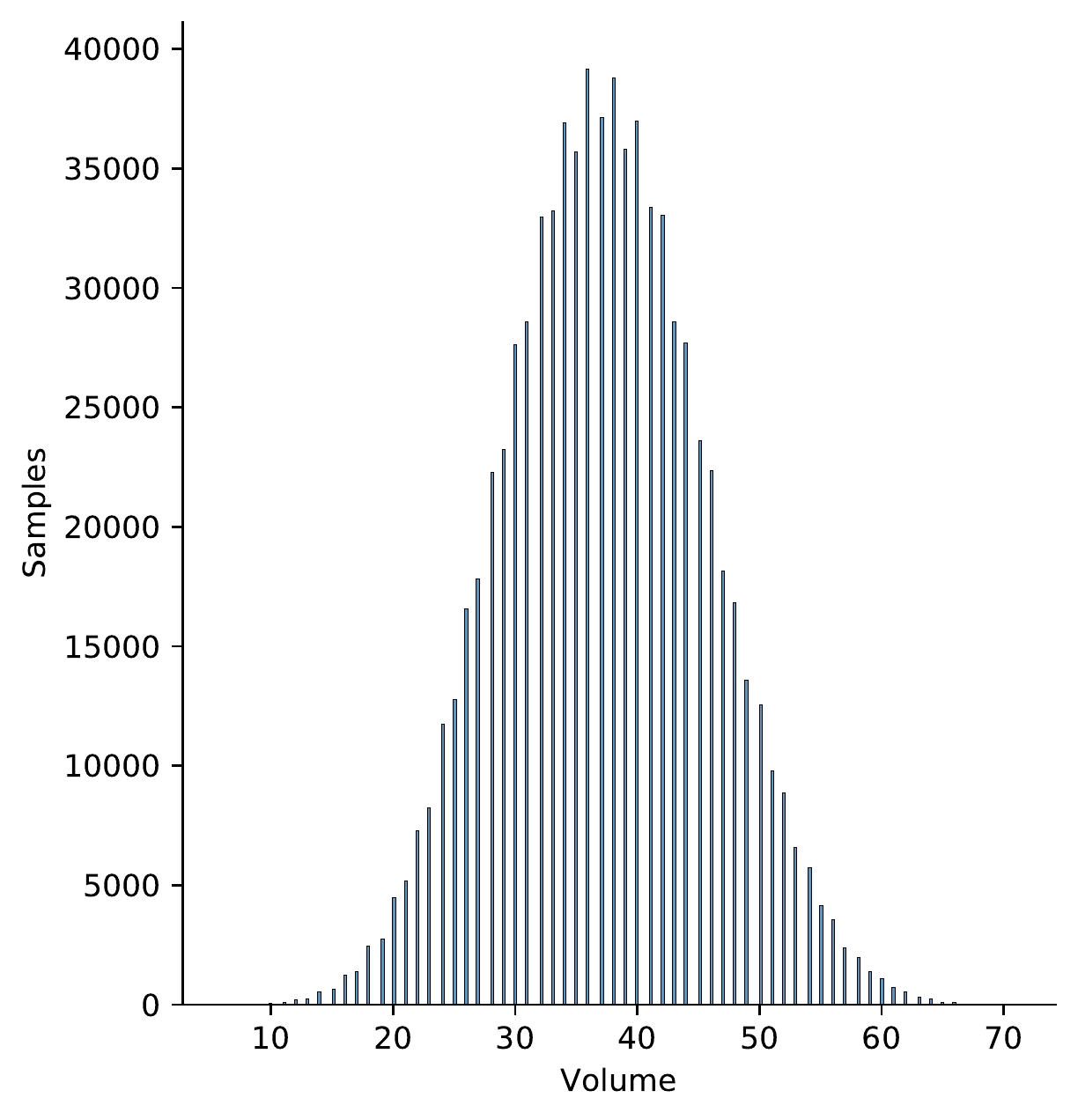}
    \caption{}\label{vol3dhist}
	\end{subfigure}
	\begin{subfigure}{6cm}
	\includegraphics[width=6cm]{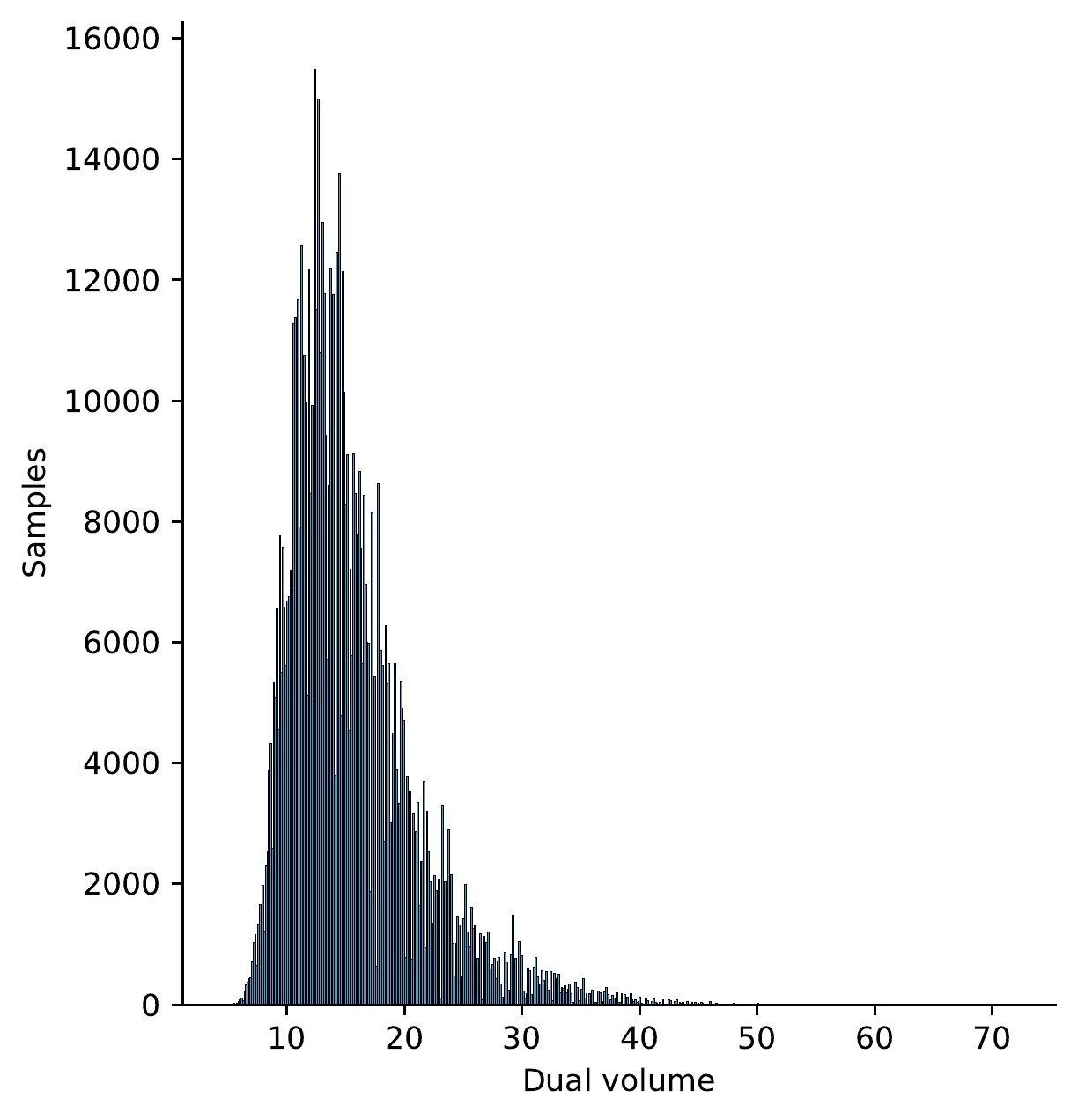}
    \caption{}\label{dualvol3dhist}
	\end{subfigure}
	\caption{(a) The distribution of the volumes for 800000 random samples. (b) The distribution of the dual volumes for 800000 random samples.}
\end{figure}

The distribution of the volumes is plotted in Figure \ref{vol3dhist}.
Similarly, we plot the distribution of the dual volumes in Figure \ref{dualvol3dhist}. Unlike the volumes which have a symmetric distribution, the dual volumes have a cluster of data between 10 and 20 while there are much fewer samples for large dual volumes.

\begin{figure}[h]
	\centering
	\begin{subfigure}[h]{0.45\textwidth}
	\centering
	\includegraphics[width=5cm]{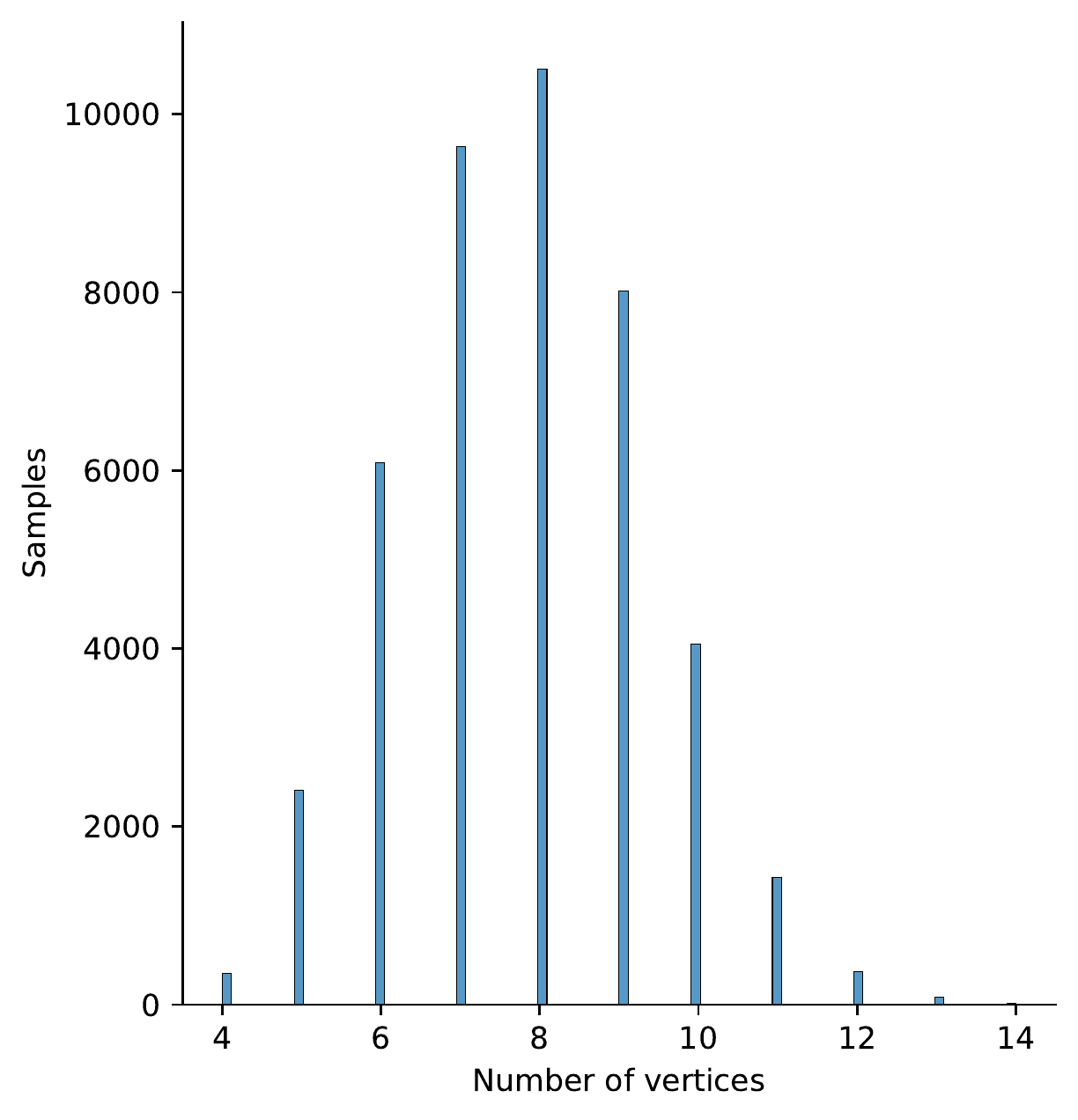}
	\caption{}
	\end{subfigure}
	\begin{subfigure}[h]{0.45\textwidth}
	\centering
	\includegraphics[width=5cm]{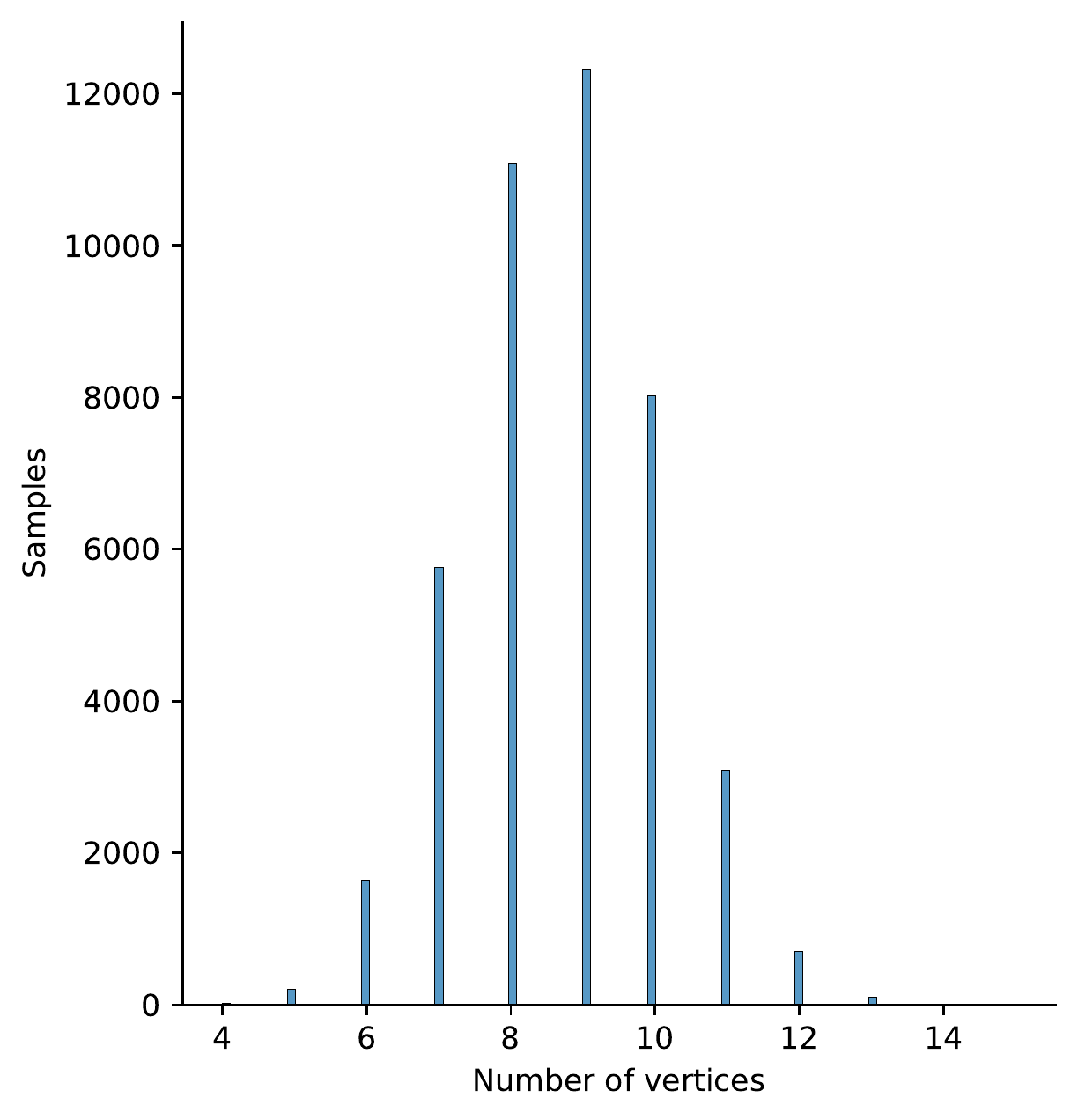}
	\caption{}
	\end{subfigure}
    \caption{(a) The distribution of the data points corresponding to reflexive polytopes, in terms of the numbers of their vertices. (b) The distribution of the 42943 randomly chosen data points corresponding to non-reflexive polytopes, in terms of the numbers of their vertices.}\label{ref3dhist}
\end{figure}

A famous subset of the canonical Fano polytopes are the 4319 reflexive polytopes \cite{Kreuzer:1998vb}, whose geometry was explored in \cite{He:2017gam}.
In our enhanced dataset, we have 42943 data points that correspond to these reflexive polytopes. For our ML task in \S\ref{ref3d}, we need to check whether the reflexive and non-reflexive ones distribute differently. Therefore, we randomly choose 42943 data points corresponding to non-reflexive polytopes and plot the distributions in terms of the numbers of vertices for the two cases in Figure \ref{ref3dhist}.
As we can see, these two plots have similar distributions.

\section{Machine Learning}\label{ml}
In this section we initiate the study of applying ML to the data introduced above.
In particular, we consider the labelled data, each entry being of the form
\begin{equation}
    \mbox{Pl\"ucker}(P) \longrightarrow
    \mbox{Property}(P)
\end{equation}
where ``property'' is one of: volume, dual volume, reflexivity, codimension, and Gorenstein index. 
This is a supervised ML problem, where we will train a NN on some percentage of data and validate on the remaining.
In this paper, we shall be examining the two lowest dimensional cases, subdividing investigations into 2d polygons\footnote{
Since there are only 16 reflexive polygons, reflexivity considerations only apply for 3d and above.
}, and 3d polytopes.

\subsection{Methodology}
\paragraph{MLP for 2d polygons} In each ML investigation the respective datasets are split into 5 subsets for a 5-fold cross-validation analysis, training 5 independent NNs each on a different complement of a fifth of the data set aside for independent validation.
The NNs trained are standard dense feed-forward, each with 4 dense layers of 64 neurons, and a leaky-ReLU activation (with $\alpha=0.01$). 
The final layer is a dense layer with a single neuron outputting the predicted property for each input vector.

They are trained in batches of 32 for 20 epochs, using an Adam optimiser \cite{kingma2017adam} to minimise a log(cosh) loss function.
In addition to the mean absolute error (MAE) learning measure for the supervised regression problem style, accuracy bins were also used to determine the ratio of test data that was within a specified range of the true value.
Note that although simple MAE was used for assessing learning performance, log(cosh) was used as the loss function in training as a smoother equivalent to improve training.

\paragraph{MLP for 3d polytopes} For 3d polytopes, we used a slightly different structure of MLP since the Pl\"ucker coordinates can be much longer. Here, we use one hidden layer with 100 perceptrons. The leaky ReLU activation function now has $\alpha=10^{-5}$. The data is trained in batches of 16 with an Adam optimizer to minimize the mean squared error loss.

\paragraph{CNN for reflexivity} Besides MLP, we also use a CNN for the classification problem of reflexivity for 3d polytopes, where we have one convolution layer followed by a Leaky ReLU layer ($\alpha=0.1$) and a MaxPooling layer. We add a SoftMax layer at the end of this structure. The data is trained in batches of 16 for 100 epochs with Adam as the optimizer to minimize the mean squared error loss.

\paragraph{Random forest for reflexivity} We can further perform the same tests using random forest, where we have 70 estimators for the classifier.


\subsection{Polygons}\label{polygon_ml}
To motivate the choice of Pl\"ucker coordinates, over the naive list of vertex coordinates, as the polygons' representation method, an initial investigation is carried out.
In this ML investigation the full dataset of polygons is considered to learn the respective polygon's volume from each input form, for each number of vertices considered. 

The NNs are trained as regressors, then assessed using varying measures.
The first measure calculates the MAE between predicted and true volumes on the test data, where a value of 0 indicates perfect learning. 
The binned accuracy measures determine the ratio of test polygons whose predicted volume is within some bin centred on the true value.
These bins are chosen to have three sizes: width 1, width $0.025 \times \text{range}$, and $0.05 \times \text{range}$; where ``range'' is the difference between the maximum and minimum volumes of polygons in the full dataset.
Accuracies of 1 indicated perfect prediction within the bin size error.
For volume, the range takes the value of 511, but is lower for other properties.

These bin styles are used across all subsequent investigations in this subsection.
The bin width of 1 identifies where properties can be predicted to the nearest integer, particularly important for those properties restricted to take integer values.
The second two bin styles based on the ranges of the properties across the dataset give a more representative and comparable measure of learning performance.

\subsubsection{Predicting Polygon Properties}
Results for these investigations, comparing inputs of vertices and Pl\"ucker coordinates, are given in Table \ref{polygon_vertpluck}. For all investigations carried out, at each number of vertices, for MAE and all accuracy bin sizes considered, learning is better using Pl\"ucker coordinates than vertices.
Interestingly in the investigations using polygons with less vertices the test set predictions are always within the range-based bins, showing that volume can be learnt very well from the polygon Pl\"ucker coordinates.
One reason we believe is that the Pl\"ucker representation preserves the rotational invariance of the polytope, so does the polytope volume. Further to this, volume is a determinant calculation from vertices which is a long computation to machine learn, whereas Pl\"ucker coordinates already encode part of this information in the computation of the minors.
Therefore in subsequent investigations we adhere to the Pl\"ucker representation for the polygons.

\begin{table}
\centering
\addtolength{\leftskip} {-2cm}
\addtolength{\rightskip}{-2cm}
\begin{tabular}{|c|c|c|c|c|c|} 
\hline
\multirow{2}{*}{\begin{tabular}[c]{@{}c@{}}Number \\ of Vertices\end{tabular}} & \multirow{2}{*}{Representation} & \multirow{2}{*}{MAE} & \multicolumn{3}{c|}{Accuracy}\\ 
\cline{4-6} & & & $\qquad\pm 0.5\qquad$ & $\pm 0.025 \times \text{range}$ & $\pm 0.05 \times \text{range}$  \\  \hline
\multirow{2}{*}{3} & Vertices & 4.941 & 0.302 & 0.891 & 0.945 \\ 
\cline{2-6} & Pl\"ucker  & 0.209  & 0.827 & 1.000 & 1.000 \\ \hline
\multirow{2}{*}{4} & Vertices  & 10.012 & 0.072 & 0.891 & 0.945 \\ 
\cline{2-6} & Pl\"ucker  & 0.625   & 1.000  & 1.000   & 1.000  \\ \hline
\multirow{2}{*}{5}  & Vertices & 8.640 & 0.060  & 0.777 & 0.926  \\ 
\cline{2-6} & Pl\"ucker & 1.051  & 0.451  & 1.000    & 1.000   \\ \hline
\multirow{2}{*}{6}  & Vertices  & 14.947   & 0.036   & 0.603  & 0.826  \\ 
\cline{2-6}  & Pl\"ucker    & 3.359   & 0.139 & 0.969 & 0.997  \\ \hline
\end{tabular}
\caption{Machine learning polygon volume from polygon representations: flattened list of vertices, or Pl\"ucker coordinates. ML runs on Dense NN Regressors, and investigations are carried out for each subset of polygons with each number of vertices. Learning is measured using MAEs, and accuracies (to 3d.p.) of test set predicted volumes being within some bin centred on the true value. The first bin has width 1, the second and third have widths of $5\% \text{ or } 10\% \times \text{range}$ respectively, where ``range'' is the difference between the maximum and minimum volumes in the full dataset (511). Pl\"ucker coordinates consistently performs better as the input representation method.}\label{polygon_vertpluck}
\end{table}

Results showing the learning accuracies (with the bin structure as described previously) for the four properties considered are given in Table \ref{polygon_ml_results}.
The ranges used in calculating the latter bin widths are: 511, 15.16, 29, 40; for volume, dual volume, Gorenstein index, and codimension respectively.
Note that the volume learning results from Pl\"ucker coordinates are repeated in both tables \ref{polygon_vertpluck} and \ref{polygon_ml_results}, for ease of comparison.
The results show that volume can be learnt especially well for all number of vertices.
In addition dual volume can be learnt well, and codimension with some success.
Conversely Gorenstein index could not be learnt with this NN architecture.

Whereas volume is a determinant and likely a simpler function of Pl\"ucker coordinates to learn, hence the good performance. 
Dual volume is less simple to compute and surprisingly still is learnt well, hinting at some simpler function that connects these inputs and outputs that shortcuts computation of the dual polytope.
The lesser performance for codimension and Gorenstein index indicate they are likely incredibly complex, if not even non-existent, functions of the Pl\"ucker coordinates (and all performed even worse using input as vertices).
In general learning is more successful for lower numbers of vertices, perhaps due to the smaller Pl\"ucker vector inputs, and the larger datasets available for these polygons.
Presumably this also explains the occasional lower performances learning with triangles compared to quadrilaterals and pentagons.

\begin{table}
\centering
\begin{tabular}{|c|c|c|c|c|c|} 
\hline
\multirow{2}{*}{Property} & \multirow{2}{*}{\begin{tabular}[c]{@{}c@{}}Number \\ of Vertices\end{tabular}} & \multirow{2}{*}{MAE} & \multicolumn{3}{c|}{Accuracy} \\ 
\cline{4-6}
 & & & $\qquad \pm 0.5\qquad$ & $\pm 0.025 \times \text{range}$ & $\pm 0.05 \times \text{range}$  \\ 
\hline
\multirow{4}{*}{Volume}   & 3  & 0.209   & 0.826 & 1.000 & 1.000 \\ 
\cline{2-6}
 & 4  & 0.615   & 0.625 & 1.000 & 1.000 \\ 
\cline{2-6}
 & 5  & 1.051   & 0.452 & 1.000 & 1.000 \\ 
\cline{2-6}
 & 6  & 3.359   & 0.139 & 0.969 & 0.997 \\ 
\hline
\multirow{4}{*}{\begin{tabular}[c]{@{}c@{}}Dual\\ Volume\end{tabular}}  & 3  & 1.181   & 0.370 & 0.370 & 0.501 \\ 
\cline{2-6}
 & 4  & 0.642   & 0.634 & 0.634 & 0.754 \\ 
\cline{2-6}
 & 5  & 0.818   & 0.496 & 0.496 & 0.638 \\ 
\cline{2-6}
 & 6  & 0.941   & 0.405 & 0.405 & 0.557 \\ 
\hline
\multirow{4}{*}{\begin{tabular}[c]{@{}c@{}}Gorenstein\\ index\end{tabular}} & 3  & 5.710   & 0.039 & 0.064 & 0.132 \\ 
\cline{2-6}
 & 4  & 5.002   & 0.069 & 0.101 & 0.196 \\ 
\cline{2-6}
 & 5  & 4.632   & 0.071 & 0.102 & 0.202 \\ 
\cline{2-6}
 & 6  & 5.343   & 0.056 & 0.080 & 0.150 \\ 
\hline
\multirow{4}{*}{Codimension} & 3  & 1.897   & 0.192 & 0.361 & 0.615 \\ 
\cline{2-6}
 & 4  & 2.726   & 0.140 & 0.268 & 0.496 \\ 
\cline{2-6}
 & 5  & 3.182   & 0.103 & 0.210 & 0.404 \\ 
\cline{2-6}
 & 6  & 2.884   & 0.126 & 0.251 & 0.468 \\
\hline
\end{tabular}
\caption{ML results for each of the properties considered, using input Pl\"ucker coordinates representing the polygons in each subset of the full dataset based on the number of polygon vertices. Learning is measured with MAE and accuracies based on test set predictions being within some bin centred on the true value. The bin widths are based on the range of values each property can take in the full dataset. Results show volume and dual volume learnt well, codimension learnt to some extent, and Gorenstein index could not be learnt. For reference the respective ranges are: 511, 15, 29, 40 for volume, dual volume, Gorenstein index, and codimension respectively.}
\label{polygon_ml_results}
\end{table}

\begin{figure}[h!]
	\centering
	\begin{subfigure}[h]{0.3\textwidth}
	\centering
	\includegraphics[width=\textwidth]{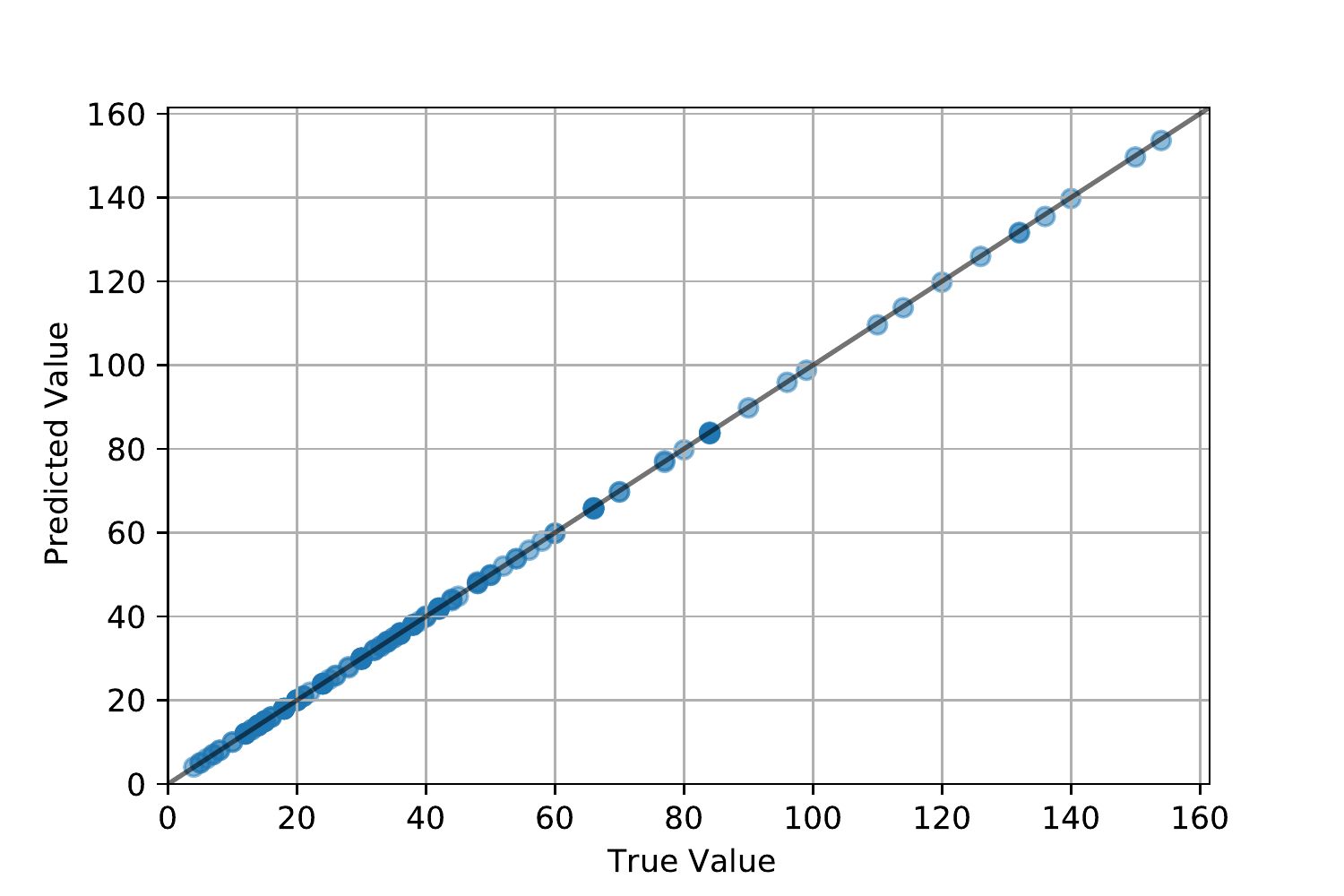}
	\caption{Volume $n=3$, MAE: 0.209}
	\end{subfigure}
	\centering
	\begin{subfigure}[h]{0.3\textwidth}
	\centering
	\includegraphics[width=\textwidth]{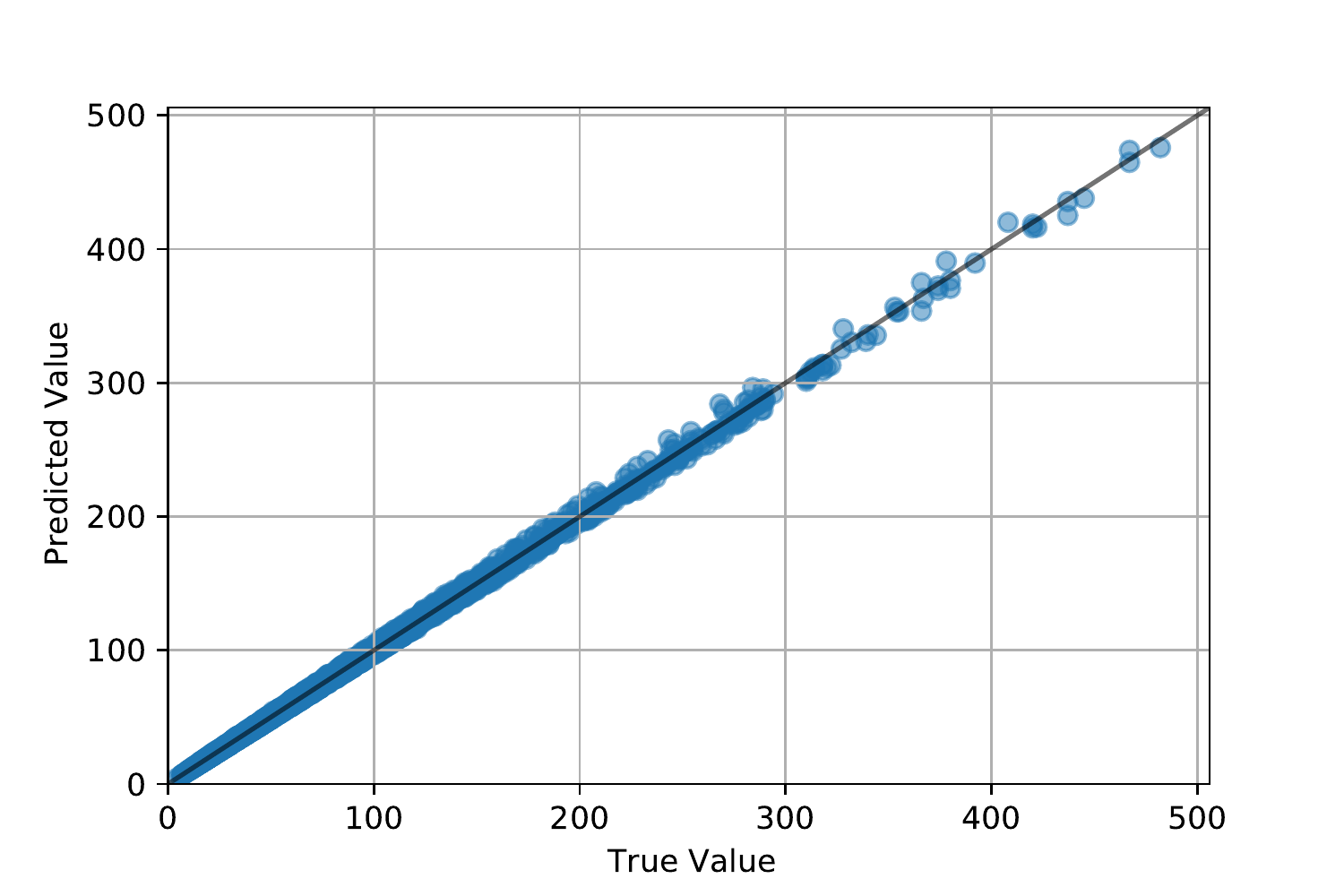}
	\caption{Volume $n=4$, MAE: 0.615}
	\end{subfigure}
	\begin{subfigure}[h]{0.3\textwidth}
	\centering
	\includegraphics[width=\textwidth]{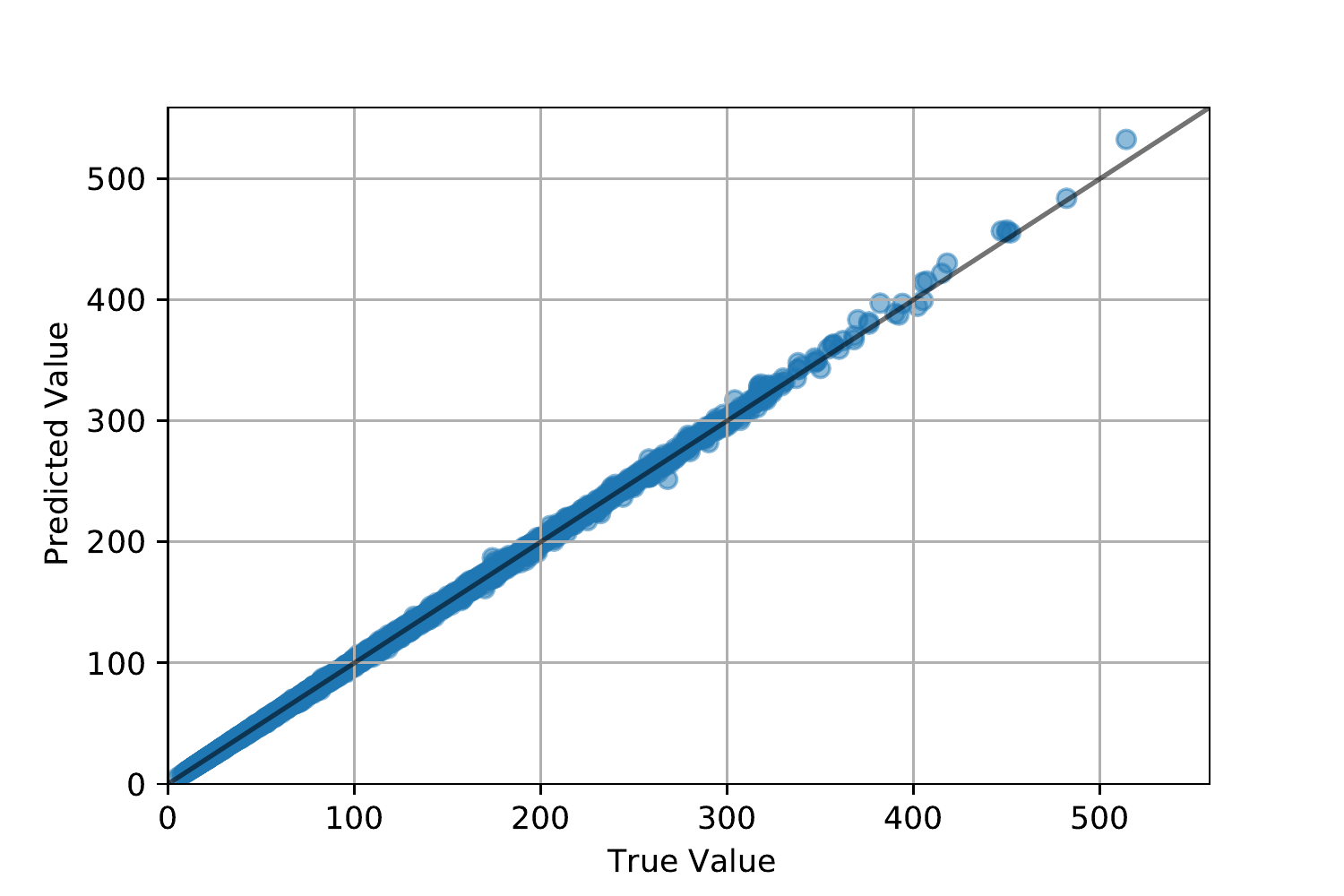}
	\caption{Volume $n=5$, MAE: 1.051}
	\end{subfigure}
	\centering
	\begin{subfigure}[h]{0.3\textwidth}
	\centering
	\includegraphics[width=\textwidth]{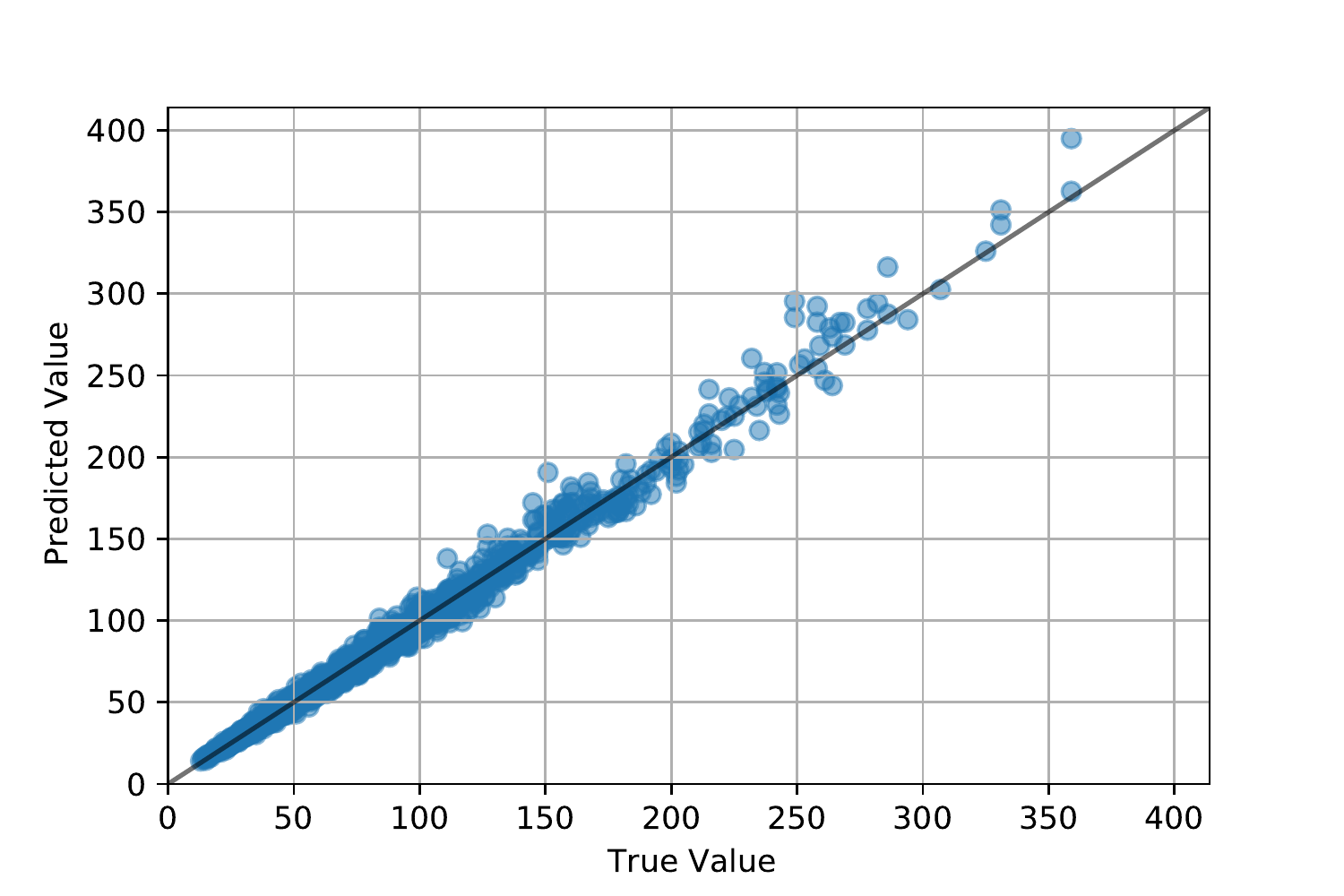}
	\caption{Volume $n=6$, MAE: 3.359}
	\end{subfigure}
	\centering
	\begin{subfigure}[h]{0.3\textwidth}
	\centering
	\includegraphics[width=\textwidth]{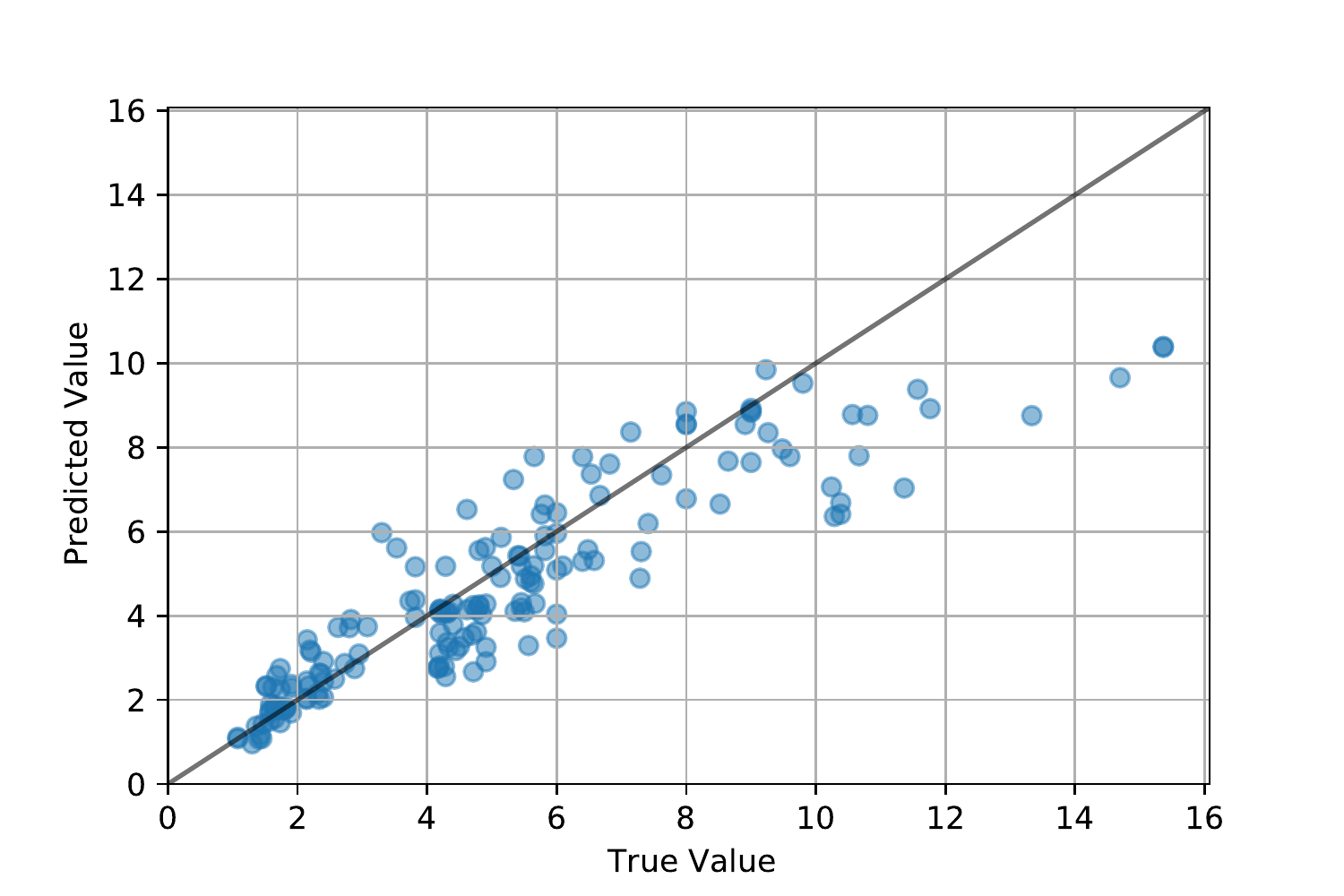}
	\caption{Dual Volume $n=3$, MAE: 1.181}
	\end{subfigure}
	\centering
	\begin{subfigure}[h]{0.3\textwidth}
	\centering
	\includegraphics[width=\textwidth]{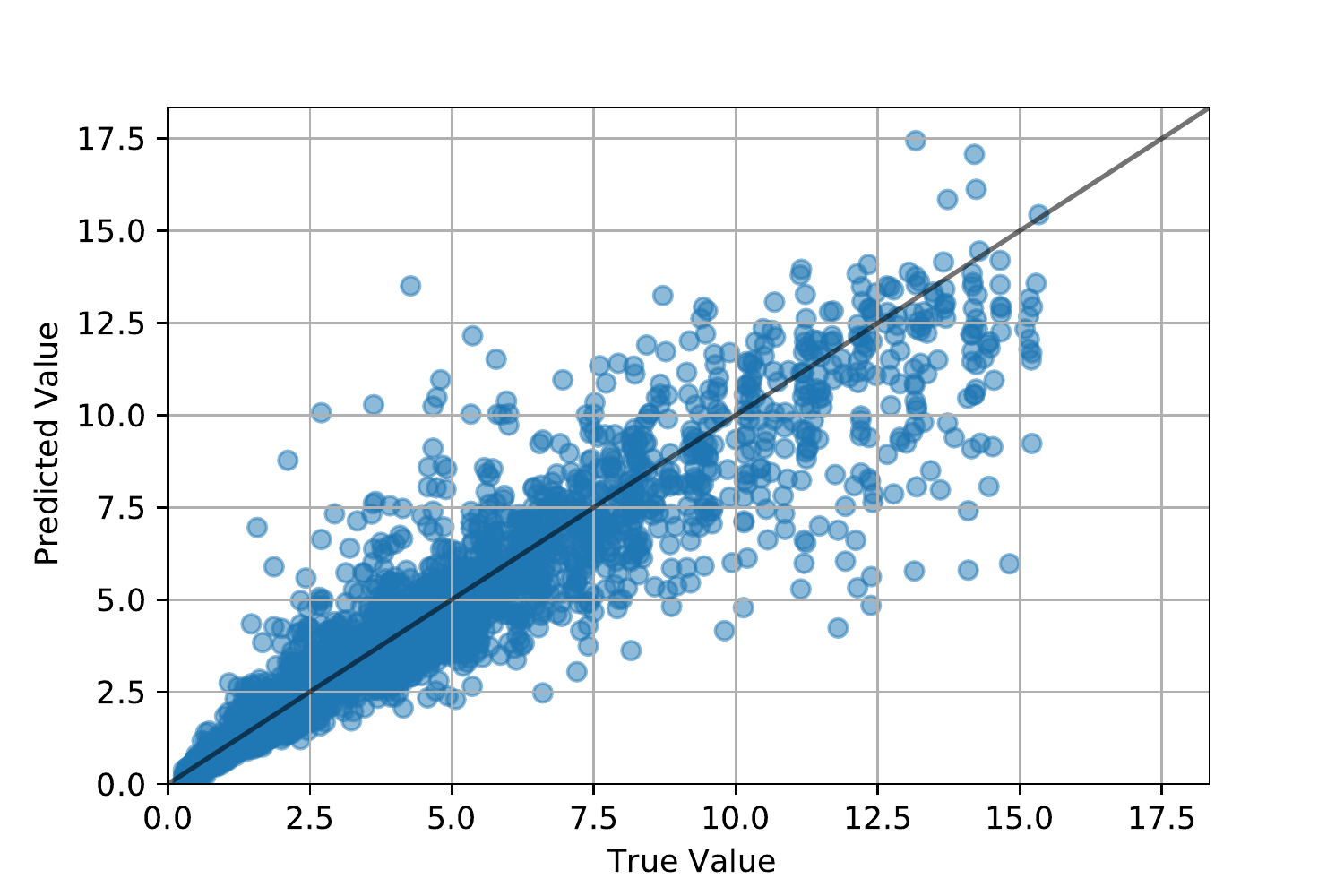}
	\caption{Dual Volume $n=4$, MAE: 0.642}
	\end{subfigure}
	\begin{subfigure}[h]{0.3\textwidth}
	\centering
	\includegraphics[width=\textwidth]{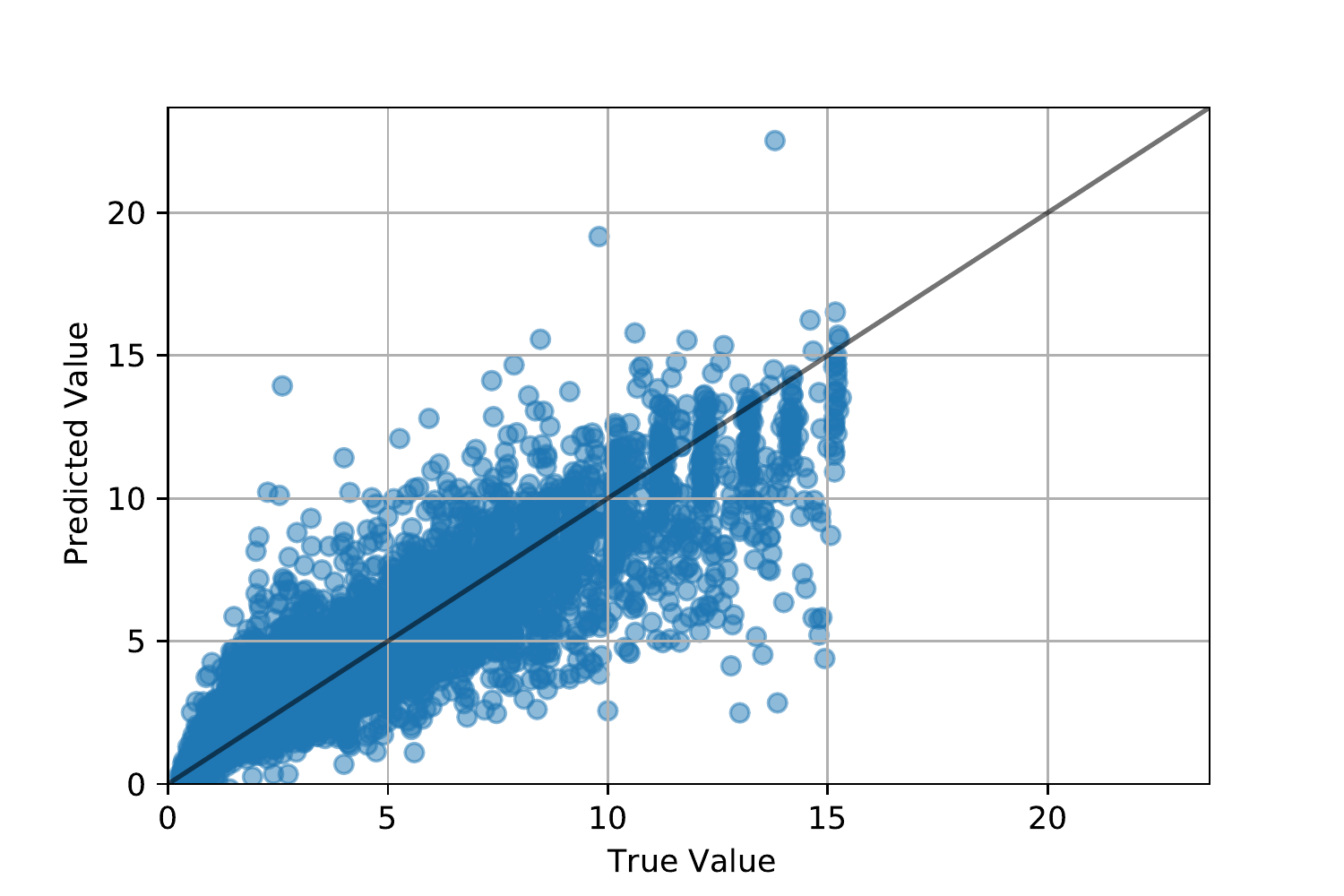}
	\caption{Dual Volume $n=5$, MAE: 0.818}
	\end{subfigure}
	\centering \quad
	\begin{subfigure}[h]{0.3\textwidth}
	\centering
	\includegraphics[width=\textwidth]{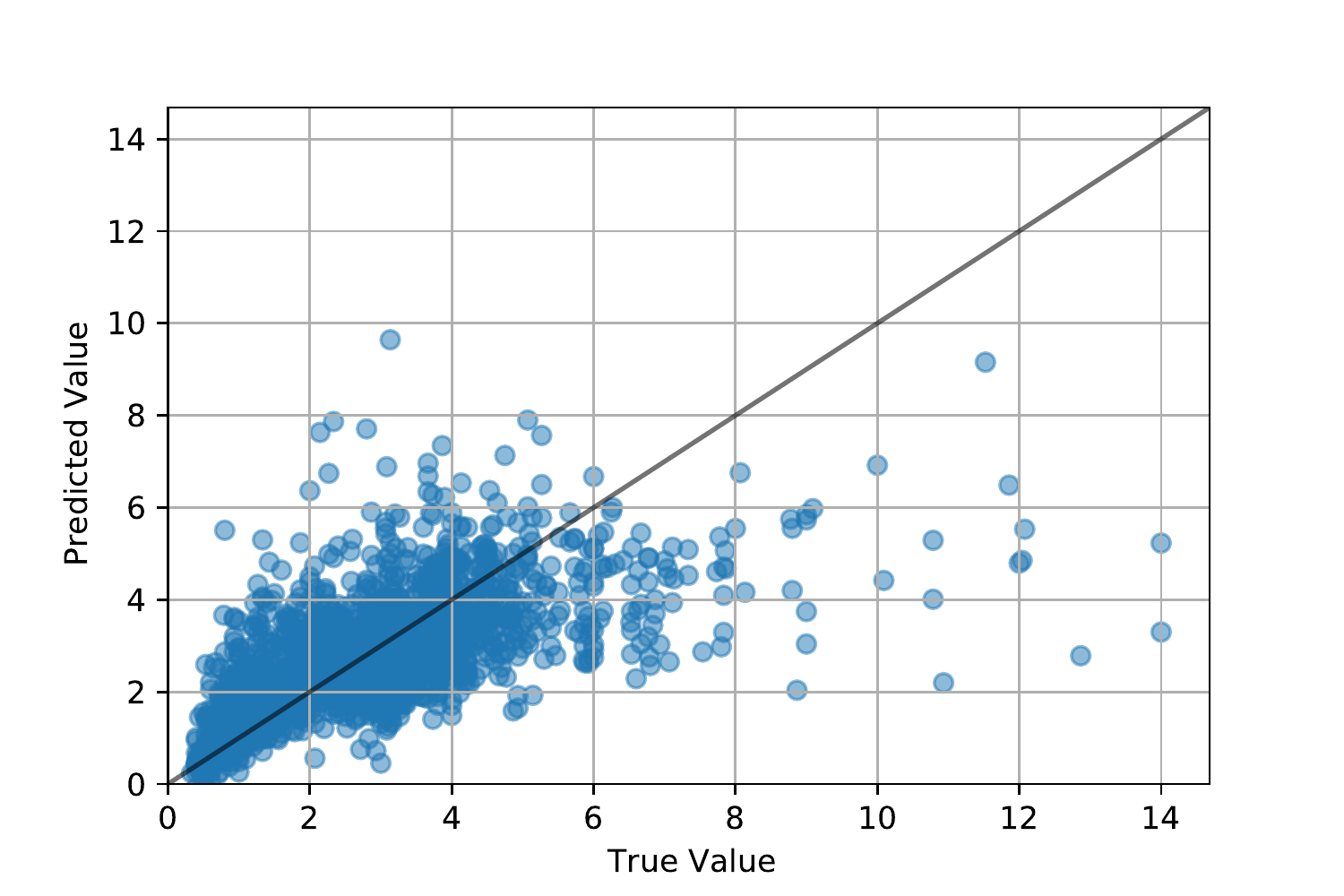}
	\caption{Dual Volume $n=6$, MAE: 0.941}
	\end{subfigure}
\caption{The true values vs NN predictions for ML volume (a-d), and dual volume (e-h), over the range of number of vertices, $n$, considered. Plots show the $y=x$ line for ease of comparison; and the MAE values are also given for each plot in the respective subcaptions to 3 decimal places.}\label{truepred2d}
\end{figure}

To further represent these learning results, the NN regressor predictions for the more successful learning of volume and dual volume are plotted against the true values over the range of numbers of vertices considered.
These plots are shown in Figure \ref{truepred2d} (plotted with the $y=x$ line for comparison), and show how a NN with the structure as before trained on 80\% of the data performed on the complement 20\% test dataset.
In addition to the plotted true vs predicted values over the respective test datasets, the respective MAEs are repeated in the respective captions of Figure \ref{truepred2d} also.

These plots reiterate the strong learning results of particularly volume, but also dual volume, from the Pl\"ucker coordinates.

\subsubsection{Pl\"ucker Augmentation}
It is anticipated the function that the NNs are trying to approximate, perhaps relies on computation of greatest common divisors (gcds).
Therefore to test this idea, the input vectors of Pl\"ucker coordinates for the pentagons (the largest dataset) are augmented with gcds in two different ways to see if the NNs could improve their learning by using these values directly.
The two gcd schemes considered used either: all $\binom{l}{2}$ pairwise gcds (denoted 2-gcds), or all gcds of the coordinate subsets formed from omitting one coordinate (denoted $(l-1)$-gcds); note $l$ is the length of the Pl\"ucker vector.

The learning results are shown in Table \ref{polygon_augmentation}, learning each of the properties with either: no augmentation, 2-gcds augmentation, or $(l-1)$-gcds augmentation.
Any learning improvement for either scheme is negligibly improved if at all, indicating that gcds were not especially important in the underlying property calculation from the input Pl\"ucker coordinates.
Note the no augmentation cases are repeated from previous tables for ease of comparison.

\begin{table}
\centering
\begin{tabular}{|c|c|c|c|c|c|} 
\hline
\multirow{2}{*}{Property} & \multirow{2}{*}{\begin{tabular}[c]{@{}c@{}}GCD\\ Scheme\end{tabular}} & \multirow{2}{*}{MAE} & \multicolumn{3}{c|}{Accuracy}   \\ 
\cline{4-6}
 & &  & $\qquad \pm 0.5 \qquad$ & $\pm 0.025 \times \text{range}$ & $\pm 0.05 \times \text{range}$  \\ 
\hline
\multirow{3}{*}{Volume}   & none & 1.051 & 0.452 & 1.000 & 1.000  \\ 
\cline{2-6}
 & $2$-gcds & 1.246 & 0.356 & 0.999  & 1.000  \\ 
\cline{2-6}
 & $(l-1)$-gcds & 0.949 & 0.455 & 1.000  & 1.000  \\ 
\hline
\multirow{3}{*}{\begin{tabular}[c]{@{}c@{}}Dual\\ Volume\end{tabular}} & none & 0.818 & 0.496 & 0.496  & 0.638  \\ 
\cline{2-6}
 & $2$-gcds & 0.865 & 0.465 & 0.465 & 0.608  \\ 
\cline{2-6}
 & $(l-1)$-gcds   & 0.832 & 0.487 & 0.487 & 0.629  \\ 
\hline
\multirow{3}{*}{\begin{tabular}[c]{@{}c@{}}Gorenstein\\ Index\end{tabular}} & none & 4.632 & 0.071    & 0.102  & 0.202  \\ 
\cline{2-6}
 & $2$-gcds  & 4.685 & 0.080 & 0.114 & 0.219  \\ 
\cline{2-6}
 & $(l-1)$-gcds & 4.489 & 0.082 & 0.119 & 0.234  \\ 
\hline
\multirow{3}{*}{Codimension}   & none & 3.182 & 0.103    & 0.210  & 0.404  \\ 
\cline{2-6}
 & $2$-gcds  & 2.569 & 0.135 & 0.263  & 0.493  \\ 
\cline{2-6}
 & $(l-1)$-gcds   & 3.215 & 0.103 & 0.205 & 0.399  \\
\hline
\end{tabular}
\caption{ML Results for learning the \textit{pentagon} properties from Pl\"ucker coordinate inputs with varying augmentations schemes based on GCDs. The learning shows the results for no augmentation, marked `none', then for augmentation with the Pl\"ucker pairwise gcds ($2$-gcds), then for augmentation with the Pl\"ucker gcds of all coordinates excluding one coordinate each time ($(l-1)$-gcds), where $l$ is the length of the Pl\"ucker coordinate vector. Learning is measured with MAE and accuracies based on the test set predictions being within some bin centred on the true value.}
\label{polygon_augmentation}
\end{table}

Further data input reformulation through one-hot encoding and principle component analysis (PCA) were also both attempted, neither leading to non-negligible improvements in the aforementioned results for learning the considered properties.

\subsubsection{Restricting Training}
\begin{figure}[h!]
	\centering
	\includegraphics[width=8cm]{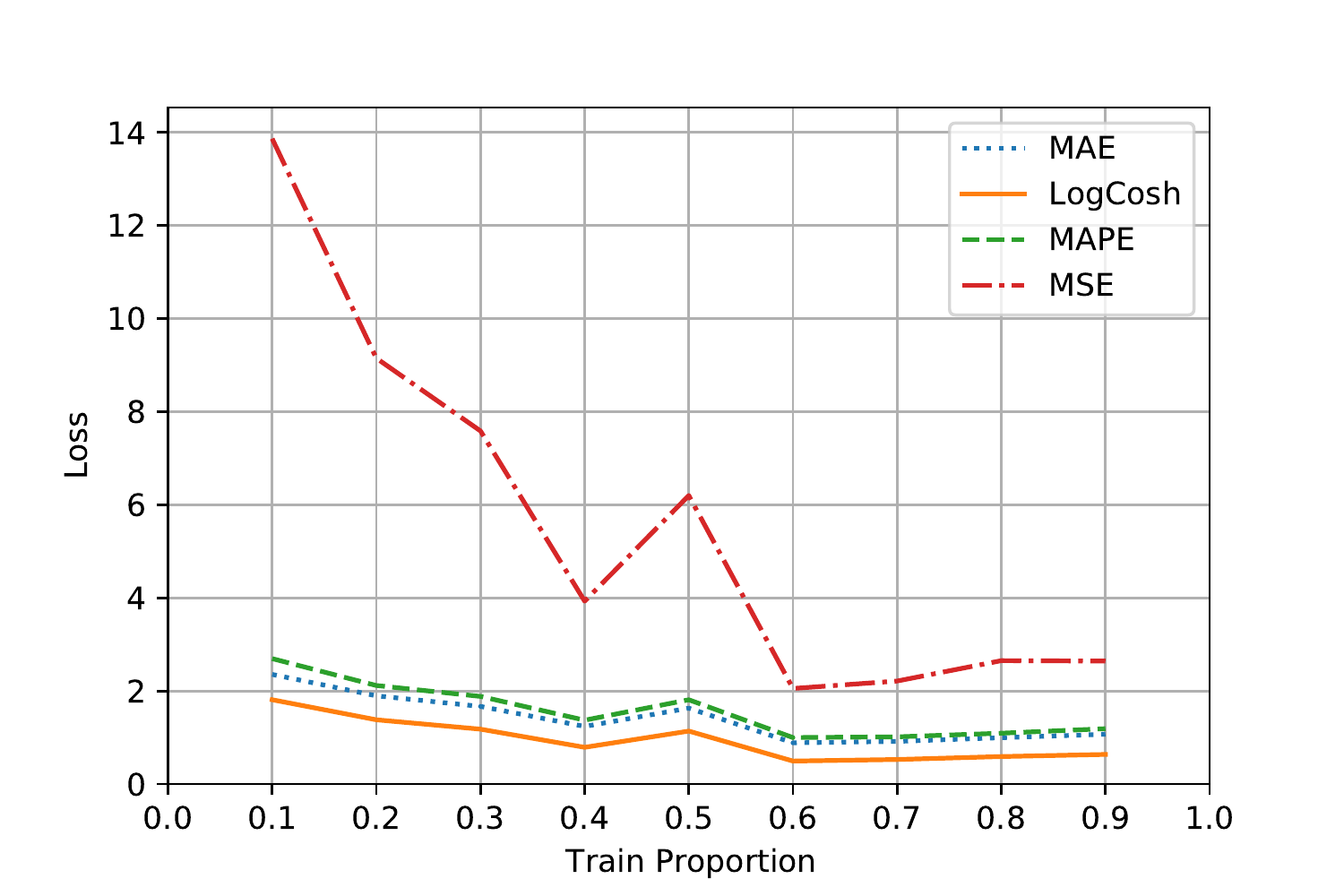}
    \caption{NN regressor measures for learning \textit{pentagon} volume from Pl\"ucker coordinates, at varying proportions of the dataset used for the training data. The results show strong learning for the standard MAE measure used such that good predictions can be achieved quickly; whilst the improvement for MSE shows that more training data does lead to better learning.}\label{pentagon_TTvary}
\end{figure}

To further examine the learning, the ML of volume from Pl\"ucker coordinates is examined at different train:test ratios. 
The results for pentagons, are shown in Figure \ref{pentagon_TTvary}, the performance at each ratio is assessed with a variety of learning measures, including: MAE, Log(cosh), Mean Absolute Percentage Error (MAPE), and Mean Squared Error (MSE). 
The NNs are still trained on the Log(cosh) loss function, which is well approximated by the more easily interpretable MAE used to assess learning in previous investigations. 
In addition the MAPE and MSE measures provide further assessments of learning performance and are plotted with those aforementioned.

Learning is still very strong when smaller datasets are used for training, since the MAE regressor measure remains low consistently as training ratio increases.
The MSE measure does however take a significantly larger value for lower accuracies, highlighting that more training data does lead to improved learning performance, as perhaps expected.
The good performance from smaller sized training datasets reinforces the idea of there being a simple function from Pl\"ucker coordinates to volume that can quickly be approximated by the NNs without much input data requires.

\subsubsection{The Inverse Problem}\label{inverse}
In vain of this success, let us consider the inverse problem of learning one of the Pl\"ucker coordinates of the polygon, given the remaining coordinates and the volume of the polygon as input data. The NN input is now thus the vector of Pl\"ucker coordinates with final coordinated omitted and the polygon volume augmented onto the end.

Motivated by previous investigations, due to the apparent existence of a closed form expression for the volume that can be inverted, we hence expect this problem to be machine-learnable. Indeed, this expectation is realised by the NN with the same architecture as before. For the polygons, we achieved MAEs $>3.5$ consistently, noting that the range in the final Pl\"ucker entry is 173 this is an impressive result. Also, accuracies are $>90\%$ for an error margin of 0.05 times the range of the respective Pl\"ucker coordinate predicted. These results are summarised in table \ref{inverse_2D}.

\begin{table}
\centering
\begin{tabular}{|c|c|c|c|c|} 
\hline
\multirow{2}{*}{Coordinates}  & \multirow{2}{*}{Number of vertices} & \multirow{2}{*}{MAE} & \multicolumn{2}{c|}{Accuracy}       \\ 
\cline{4-5}
&       &  & $\pm 0.025 \times \text{range}$ & $\pm 0.05 \times \text{range}$  \\ 
\hline
\multirow{4}{*}{\begin{tabular}[c]{@{}c@{}}Pl\"ucker\\ coordinates\end{tabular}} & 3     & 0.295      & 1.000       & 1.000       \\ 
\cline{2-5}
& 4     & 2.599      & 0.950       & 0.983       \\ 
\cline{2-5}
& 5     & 2.061      & 0.981       & 0.994       \\ 
\cline{2-5}
& 6     & 3.169      & 0.900       & 0.971       \\ 
\hline
\end{tabular}
\caption{ML Results for learning one of the Pl\"ucker coordinates from the rest of the Pl\"ucker coordinates and the volumes of the polygons. Results show MAE (note the predicted final Pl\"ucker coordinate has a range of 173), and the binned accuracies.}
\label{inverse_2D}
\end{table}

\subsection{Polytopes}\label{plkpolyhedra}
Following success learning polygon properties in 2d, we can turn our attention to performing similar tests on polytopes in 3d as well. The input vectors are still the Pl\"ucker coordinates, and we pad 0's to the end of the vectors to cope with varying input length. We start with learning their volumes.

\subsubsection{Predicting Polytope Volumes}\label{vol3d}
We use MLP regressors to predict the polytope volumes $V$ as $4\leq V\leq72$, now across all values for the number of vertices. We randomly choose 800000 samples from the dataset for training (whose distribution was shown in Figure \ref{vol3dhist}). At $90\%$ training percentage, the results are listed in Table \ref{plk_vol3d}.
\begin{table}[h]
\centering
\begin{tabular}{|c|c|c|c|c|c|}
\hline
MAE & \multicolumn{5}{c|}{1.680} \\ \hline\hline
Error & $\pm0.5$ & $\pm1$ & $\pm2$ & $\pm3$ & $\pm4$ \\
Accuracy & 0.196 & 0.380 & 0.667 & 0.846 & 0.936\\
\hline
\end{tabular}
\caption{The MAE and accuracies for predicting the volumes with 800000 samples at 90\% training percentage. The errors indicate the accuracies/percentages of the predictions $V_\text{pred}\in[V_\text{actual}-\text{err},V_\text{actual}+\text{err}]$.}\label{plk_vol3d}
\end{table}

As we can see, both MLP and CNN can give pretty high accuracies with $\pm5.797\%(=\pm4/69)$ error tolerance. In particular, MLP easily gives over $90\%$ accuracy. Moreover, given such a large range of volumes, the MAE result shows that the mean difference between a predicted volume and a true volume is even less than 2. One may wonder how these trained models would perform for the remaining samples among the 6.7 million data points. In fact, it is natural to expect that the models would keep their efficacy and have the same performance as in Table \ref{plk_vol3d}. This is because random samples are chosen from a fixed dataset and hence there is no extrapolation for the machine. More importantly, the results would not change since we are dealing with mathematical objects without noise in the data. If the machine finds any patterns in the training data (which is chosen randomly), all such data should also follow the patterns. Indeed, this is verified in Table \ref{plk_vol3dmore} with 100000 more samples from the remaining data\footnote{Henceforth, we will not repeat this for future tests. Unless specified, we will only list the validation results for part of the data (due to the huge amount of the data) and the efficacy of the models on \emph{any} polytope should be guaranteed following the argument here.}.
\begin{table}[h]
\centering
\begin{tabular}{|c|c|c|c|c|c|}
\hline
MAE & \multicolumn{5}{c|}{1.644} \\ \hline\hline
Error & $\pm0.5$ & $\pm1$ & $\pm2$ & $\pm3$ & $\pm4$ \\
Accuracy & 0.202 & 0.382 & 0.678 & 0.852 & 0.942\\
\hline
\end{tabular}
\caption{The MAE and accuracies for predicting the volumes with 100000 more samples.}\label{plk_vol3dmore}
\end{table}

\begin{figure}[h]
	\centering
	\includegraphics[width=7cm]{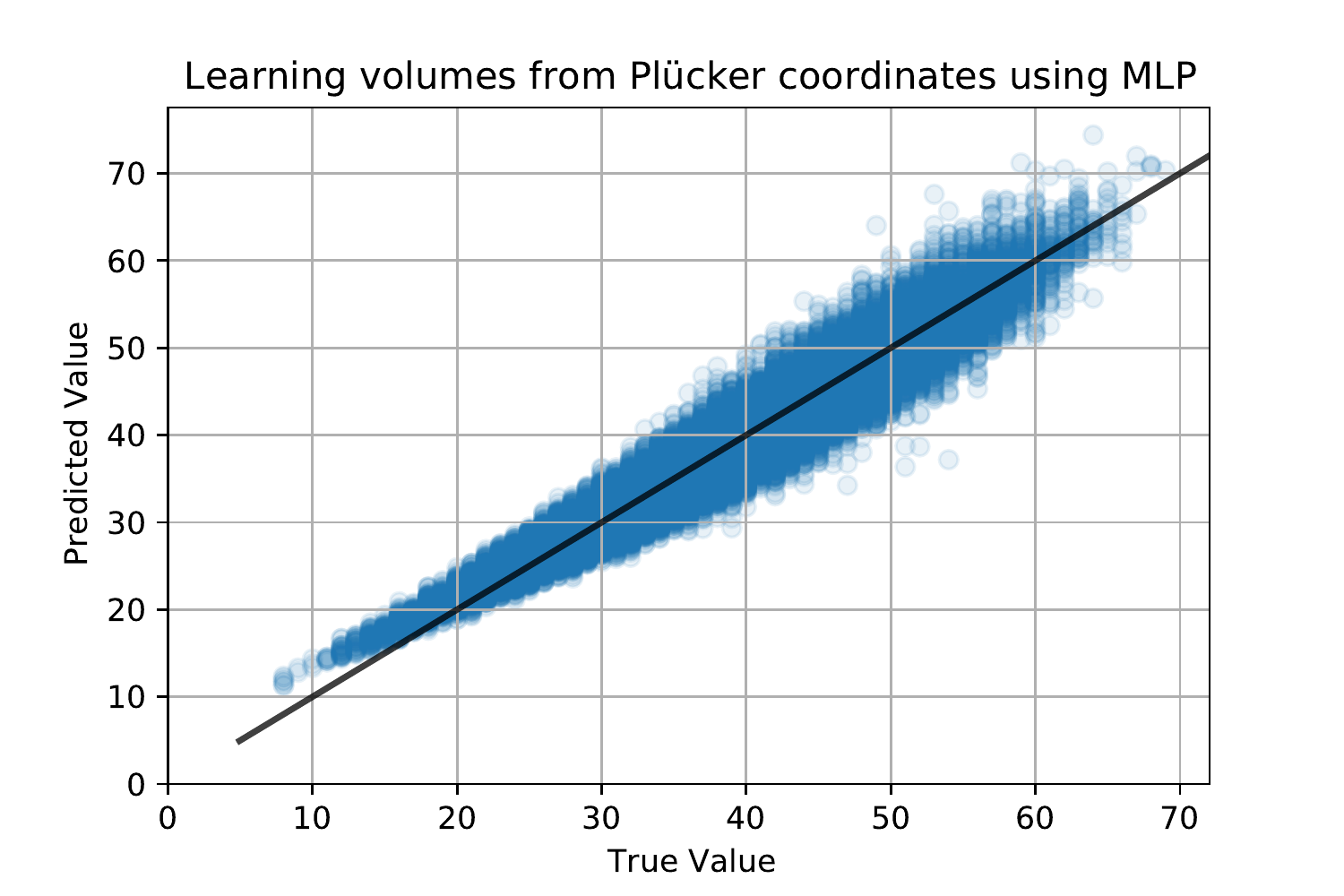}
    \caption{The true volumes vs predicted volumes. The PMCC is 0.967.}\label{truepredvol3d}
\end{figure}

We also plot the true values versus predicted values for the 80000 test samples as 10\% of the dataset in Figure \ref{truepredvol3d}.

This shows that by training only $720000/6744246=10.676\%$ of the data, we can already reach rather high accuracy for the polytope volumes.

Similar to the polygon tests before, we may also augment the input Pl\"ucker coordinates with for instance their gcds to see if the results would be improved. It turns out that including all the ``($l-1$)-gcds'' (viz, gcds of $(l-1)$ elements in an input vetor of length $l$) would give roughly the same performance as in Table \ref{gcdlen-1_vol3d}.
\begin{table}[h]
\centering
\begin{tabular}{|c|c|c|c|c|c|}
\hline
MAE & \multicolumn{5}{c|}{1.696} \\ \hline\hline
Error & $\pm0.5$ & $\pm1$ & $\pm2$ & $\pm3$ & $\pm4$ \\
Accuracy & 0.190 & 0.371 & 0.661 & 0.844 & 0.936\\
\hline
\end{tabular}
\caption{The MAE and accuracies for predicting the volumes with 720000 training samples. The input Pl\"ucker coordinates are augmented with their ``($l-1$)-gcds''.}\label{gcdlen-1_vol3d}
\end{table}
One may also try other manipulations of the input (such as augmentation with pairwise gcds or using one-hot encoding or PCAs), but they would all lead to worse results. Therefore, taking the computation cost into account, bare Pl\"ucker coordinates (without further modifications) is the optimal choice for predicting polytope volumes.

\paragraph{Pl\"ucker coordinates with fixed length} Recall that we pad 0's to the input vectors due to their different lengths. Although the Pl\"ucker coordinates may already have 0 as an entry indicating certain relations among the vertices, 0 is yet the most natural choice for padding\footnote{One may also try some other values for padding such as sufficiently large positive or negative numbers. However, all would make the predictions much worse.}. On the other hand, we may also restrict the input Pl\"ucker coordinates to some fixed length as in \S\ref{polygon_ml} so as to see whether such padding would bring noise to the input and whether fixed length in lieu of varying lengths would help the machine to learn the volumes.

We first choose Pl\"ucker coordinates of length 10. Now there are only 36487 samples of length 10 in the dataset and $6\leq V\leq70$. We train 4500 samples (i.e., $12\%$ of 36487 samples) as in Table \ref{plk10_vol3d}.
\begin{table}[h]
\centering
\begin{tabular}{|c|c|c|c|c|c|}
\hline
MAE & \multicolumn{5}{c|}{0.305} \\ \hline\hline
Error & $\pm0.5$ & $\pm1$ & $\pm2$ & $\pm3$ & $\pm4$ \\
Accuracy & 0.820 & 0.974 & 1.000 & 1.000 & 1.000\\
\hline
\end{tabular}
\caption{The MAE and accuracies for predicting the volumes with 4500 training samples. The input Pl\"ucker coordinates have length fixed to 10.}\label{plk10_vol3d}
\end{table}
We also test this with augmenting ($l-1$)-gcds in Table \ref{plk10gcdlen-1_vol3d}.
\begin{table}[h]
\centering
\begin{tabular}{|c|c|c|c|c|c|}
\hline
MAE & \multicolumn{5}{c|}{0.311} \\ \hline\hline
Error & $\pm0.5$ & $\pm1$ & $\pm2$ & $\pm3$ & $\pm4$ \\
Accuracy & 0.786 & 0.978 & 0.996 & 0.998 & 1.000\\
\hline
\end{tabular}
\caption{The accuracies for predicting the volumes with 4500 training samples. The input Pl\"ucker coordinates have length fixed to 10 and are further augmented with ($l-1$)-gcds (i.e., 9-gcds).}\label{plk10gcdlen-1_vol3d}
\end{table}

As a result, we see that the performance is greatly improved for Pl\"ucker coordinates of just length 10 compared to those for various lengths.

We may also try Pl\"ucker coordinates with longer fixed lengths. For instance, the results for those of lengths 35 and 56 are listed in Table \ref{plk35/56_vol3d}.
\begin{table}[h]
\centering
\begin{tabular}{|c|c|c|c|c|c|c|}
\hline
\multirow{2}{*}{Length 35} & MAE & \multicolumn{5}{c|}{0.908} \\ \cline{2-7}
& Accuracies & 0.351 & 0.636 & 0.916 & 0.986 & 0.997 \\ \hline
\multirow{2}{*}{Length 35 with ($l-1$)-gcds} & MAE & \multicolumn{5}{c|}{1.090} \\ \cline{2-7}
& MLP & 0.290 & 0.546 & 0.859 & 0.967 & 0.993 \\ \hline
\multirow{2}{*}{Length 56} & MAE & \multicolumn{5}{c|}{1.408} \\ \cline{2-7}
& MLP & 0.226 & 0.434 & 0.749 & 0.910 & 0.973 \\ \hline
\multirow{2}{*}{Length 56 with ($l-1$)-gcds} & MAE & \multicolumn{5}{c|}{1.543} \\ \cline{2-7}
& MLP & 0.209 & 0.404 & 0.706 & 0.877 & 0.956\\
\hline
\end{tabular}
\caption{The MAE and accuracies for predicting the volumes with input vectors of fixed lengths. The accuracies are still for error tolerance $(\pm0.5,\pm1,\pm2,\pm3,\pm4)$. For length 35, there are 899564 samples and $10\leq V\leq69$. We train 90000 samples (i.e., $10\%$ of all the samples) for the models. For length 56, there are 1728583 samples and $12\leq V\leq70$. We train 90000 samples (i.e., $5\%$ of all the samples) for the models.}\label{plk35/56_vol3d}
\end{table}
We find that the performance also has a significant improvement though shorter fixed lengths would have greater improvements. It is also worth noting that the training percentages of all the samples are similar to the percentage for the case with varying lengths. Remarkably, for the tests with fixed length 56, performance is good with the training percentage only $5\%$.

\subsubsection{Predicting Dual Volumes}\label{dualvol3d}
In addition to examining volume, the Pl\"ucker coordinates are used in this investigation to predict the dual volumes $V_\text{dual}$ of the polytopes. Again, we use MLP and CNN regressors as $4\leq V_\text{dual}\leq72$. We still train 720000 samples, with results are listed in Table \ref{plk_dualvol3d}.
\begin{table}[h]
\centering
\begin{tabular}{|c|c|c|c|c|c|c|}
\hline
MAE & \multicolumn{6}{c|}{2.590} \\ \hline\hline
Error & $\pm0.5$ & $\pm1$ & $\pm2$ & $\pm3$ & $\pm4$ & $\pm5$ \\
Accuracy & 0.129 & 0.253 & 0.478 & 0.669 & 0.807 & 0.890 \\
\hline
\end{tabular}
\caption{The MAE and accuracies for predicting the dual volumes with 720000 samples.}\label{plk_dualvol3d}
\end{table}
As we can see, the performance is not as good as for volumes. Nevertheless, the mean difference between predicted values and true values is only slightly above 2.5. The accuracy could still approach almost $90\%$ if we allow an error tolerance of $\pm7.2\%(=\pm5/69)$.

\begin{figure}[h]
	\centering
	\includegraphics[width=7cm]{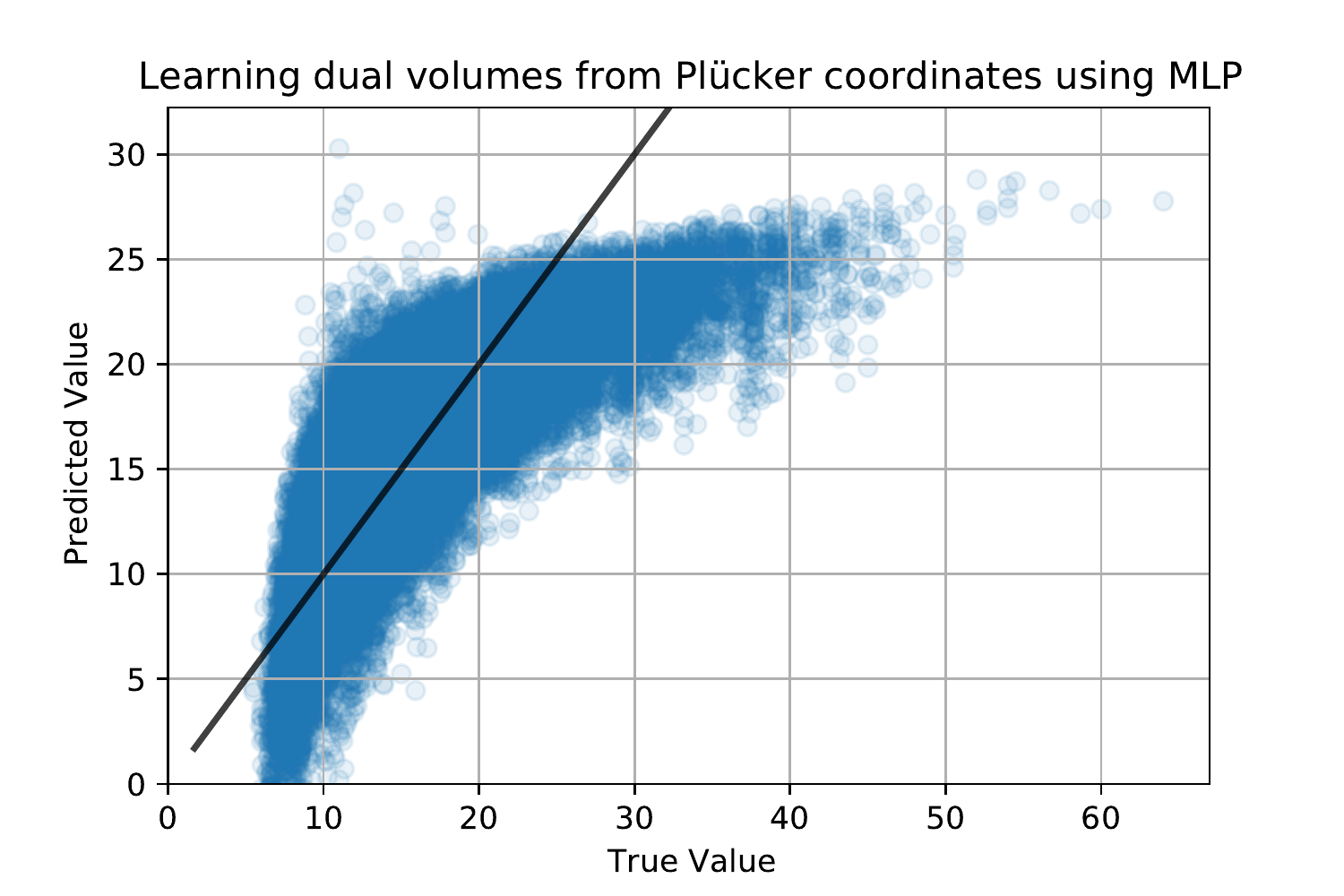}
    \caption{The true dual volumes vs predicted dual volumes. The PMCC is 0.770.}\label{truepreddualvol3d}
\end{figure}

We also plot the true values versus predicted values for 80000 random test samples in Figure \ref{truepreddualvol3d}.

We can also augment the input Pl\"ucker coordinates with their $(l-1)$-gcds as before to see how this would affect the performance. The results can be found in Table \ref{gcdlen-1_dualvol3d}.
\begin{table}[h]
\centering
\begin{tabular}{|c|c|c|c|c|c|c|}
\hline
MAE & \multicolumn{6}{c|}{2.414} \\ \hline\hline
Error & $\pm0.5$ & $\pm1$ & $\pm2$ & $\pm3$ & $\pm4$ & $\pm5$ \\
Accuracy & 0.145 & 0.280 & 0.529 & 0.713 & 0.833 & 0.902 \\
\hline
\end{tabular}
\caption{The accuracies for predicting the dual volumes with 720000 samples. The input Pl\"ucker coordinates are augmented with their $(l-1)$-gcds.}\label{gcdlen-1_dualvol3d}
\end{table}
One may also modify the input with other methods such as one-hot encoding or PCA, but they would not improve the performance as seen in previous investigations.

To check whether the data distribution has affected the above machine-learning results, recall the uneven distribution for dual volumes in Figure \ref{dualvol3dhist}. It is possible that the machine is just predicting values in the range $[10,20]$ all the time to get a near 0.9 accuracy (for $\text{err}=\pm5$). To avoid potential ``bias'' caused by this distribution, we further choose 10000 random samples whose dual volumes are greater than or equal to 25 and use the above trained models to predict their dual volumes. The results are listed in Table \ref{inputslarge_dualvol3d}.
\begin{table}[h]
\centering
\begin{tabular}{|c|c|c|c|c|c|c|c|}
\hline
\multirow{2}{*}{Pl\"ucker coordinates} & MAE & \multicolumn{6}{c|}{2.647} \\ \cline{2-8}
& Accuracy & 0.129 & 0.253 & 0.478 & 0.663 & 0.801 & 0.881 \\ \hline
\multirow{2}{*}{Augmentation with ($l-1$)-gcds} & MAE & \multicolumn{6}{c|}{2.447} \\ \cline{2-8}
& Accuracy & 0.140 & 0.274 & 0.519 & 0.708 & 0.824 & 0.900 \\
\hline
\end{tabular}
\caption{The MAEs and accuracies for predicting the dual volumes in the range $[25,72]$. The accuracies are still for error tolerance $(\pm0.5,\pm1,\pm2,\pm3,\pm4,\pm5)$.}\label{inputslarge_dualvol3d}
\end{table}
Compared to the results in Table \ref{plk_dualvol3d} and \ref{gcdlen-1_dualvol3d}, we find that the models still have quite similar performance.

\paragraph{Pl\"ucker coordinates with fixed length} As before, let us also check whether fixing the input length would improve the predictions for dual volumes. Again, we choose Pl\"ucker coordinates of length 10 as an example, where the dual volume ranges from 4.8 to 66.5. Altogether there are 36487 samples. It turns out that even at $90\%$ training as in Table \ref{inputs10_dualvol3d}, the results are not improved (even worse) compared to varying input lengths.
\begin{table}[h]
\centering
\begin{tabular}{|c|c|c|c|c|c|c|c|}
\hline
\multirow{2}{*}{Pl\"ucker coordinates} & MAE & \multicolumn{6}{c|}{3.675} \\ \cline{2-8}
& Accuracy & 0.086 & 0.171 & 0.347 & 0.501 & 0.650 & 0.752 \\ \hline
\multirow{2}{*}{Augmentation with ($l-1$)-gcds} & MAE & \multicolumn{6}{c|}{4.051} \\ \cline{2-8}
& Accuracy & 0.073 & 0.146 & 0.295 & 0.443 & 0.583 & 0.710 \\
\hline
\end{tabular}
\caption{The MAEs and accuracies for predicting the dual volumes whose Pl\"ucker coordinates are all of length 10 using different types of input vectors. The accuracies are still for error tolerance $(\pm0.5,\pm1,\pm2,\pm3,\pm4,\pm5)$.}\label{inputs10_dualvol3d}
\end{table}

Even though this is already $90\%$ of the whole data, the size of the training set is still not comparable to 72000. Therefore, it is natural to expect the worse predictions here. However, even if we had sufficient data, there would still be a failure of possible improvements. To see this, we also train the same amount of data (i.e., 32838 samples) whose Pl\"ucker coordinates have different lengths. The results are listed in Table \ref{inputscompare_dualvol3d}. As we can see, the volumes and dual volumes behave very differently in terms of Pl\"ucker coordinates.
\begin{table}[h]
\centering
\begin{tabular}{|c|c|c|c|c|c|c|c|}
\hline
\multirow{2}{*}{Pl\"ucker coordinates} & MAE & \multicolumn{6}{c|}{3.331} \\ \cline{2-8}
& Accuracy & 0.106 & 0.218 & 0.400 & 0.559 & 0.696 & 0.787 \\ \hline
\multirow{2}{*}{Augmentation with ($l-1$)-gcds} & MAE & \multicolumn{6}{c|}{3.190} \\ \cline{2-8}
& Accuracy & 0.106 & 0.208 & 0.413 & 0.582 & 0.703 & 0.803 \\
\hline
\end{tabular}
\caption{The MAEs and accuracies for predicting the dual volumes whose Pl\"ucker coordinates are of different lengths using different types of input vectors. The accuracies are still for error tolerance $(\pm0.5,\pm1,\pm2,\pm3,\pm4,\pm5)$. The training set has $32838\,(=36487\times90\%)$ data points.}\label{inputscompare_dualvol3d}
\end{table}

\subsubsection{Reflexivity of the Polytopes}\label{ref3d}
Unlike 2d lattice polygons among which only 16 are reflexive, there are 4319 reflexive polytopes in our dataset of 3d canonical polytopes. Therefore, we may also use the Pl\"ucker coordinates as input to tell whether a polytope is reflexive or not. Now this becomes a binary classification where we use 0 to denote non-reflexive polytopes and 1 for reflexive ones.

In our dataset with 6744246 polytopes, 42943 of them are reflexive. Therefore, we further choose 42943 non-reflexive ones to create a balanced dataset. We use $90\%$ of them to train the models\footnote{Recall that this is essentially training $\sim77000$ samples out of the 6.7 million data points and the results should be valid for the whole dataset as we are machine-learning mathematical objects.}. Balanced datasets prevent the machine easily reaching over 0.9 accuracy by simply predicting 0 all the time, which would happen if the full unbalanced (therefore biased) dataset is used.

Here, we compare different types of inputs including the bare Pl\"ucker coordinates, augmentation with $(l-1)$-gcds and one-hot encoding. Besides MLP and CNN classifiers, we also employ a random forest classifier here.
\begin{table}[h]
\centering
\begin{tabular}{|c|c|c|c|c|}
\hline
Classifier & Random Forest & MLP & CNN \\ \hline
Pl\"ucker coordinates & 0.764 & 0.706 & 0.746 \\ \hline
Augmentation with ($l-1$)-gcds & 0.767 & 0.734 & 0.735 \\ \hline
One-hot encoding & 0.770 & 0.790 & 0.813\\
\hline
\end{tabular}
\caption{The accuracies for classifying reflexivity using different types of input vectors.}\label{inputs_ref}
\end{table}
As we can see from Table \ref{inputscompare_ref}, one-hot encoding can give the best performance to classify reflexivity, especially when using CNN.

\paragraph{Pl\"ucker coordinates with fixed length} Again let us restrict the lengths of the Pl\"ucker coordinates to be 10 and check whether the performance of the classification can be improved. There are 2410 reflexive samples and we randomly choose same number non-reflexive ones to make the dataset balanced. The performance of the models at $90\%$ training can be found in Table \ref{inputs10_ref}.
\begin{table}[h]
\centering
\begin{tabular}{|c|c|c|c|c|}
\hline
Classifier & Random Forest & MLP & CNN \\ \hline
Pl\"ucker coordinates & 0.712 & 0.724 & 0.718 \\ \hline
Augmentation with ($l-1$)-gcds & 0.761 & 0.761 & 0.747 \\ \hline
One-hot encoding & 0.707 & 0.713 & 0.724\\
\hline
\end{tabular}
\caption{The accuracies for classifying reflexivity using different types of input vectors. The length of Pl\"ucker coordinates is fixed to be 10.}\label{inputs10_ref}
\end{table}

Again, let us also check the results with (2410+2410) samples of varying input lengths for comparison. The accuracies are shown in Table \ref{inputscompare_ref}.
\begin{table}[h]
\centering
\begin{tabular}{|c|c|c|c|c|}
\hline
Classifier & Random Forest & MLP & CNN \\ \hline
Pl\"ucker coordinates & 0.763 & 0.601 & 0.645 \\ \hline
Augmentation with ($l-1$)-gcds & 0.746 & 0.674 & 0.668 \\ \hline
One-hot encoding & 0.728 & 0.693 & 0.753\\
\hline
\end{tabular}
\caption{The accuracies for classifying reflexivity using different input vectors whose corresponding Pl\"ucker coordinates have varying lengths. There are $4338\,(=4820 \times 90\%)$ training samples.}\label{inputscompare_ref}
\end{table}
As we can see, although MLP and CNN have better performance for fixed input length, the best results amongst all the models, which are from random forest, remain roughly the same. Therefore padding 0's to the input vectors with varying length does not seem to affect the learning results for reflexivity.

\subsubsection{From Polygons to Polytopes}\label{2dto3d}
So far, we have conducted various tests for both 2d and 3d polytopes. It would be an interesting question to ask if a trained model for 2d polygons can learn the same quantities for 3d polytopes with less training data compared to a model trained with exclusively 3d data.

Let us start with the volumes. We first consider 2d and 3d polytopes having Pl\"ucker coordinates of same length. A good choice would be length equal to 10 because both 2d and 3d polytopes happen to have 5 vertices in this situation as $\binom{5}{3}=\binom{5}{2}=10$ respectively. We first train our model only with 12000 2d data points, whose good performance is expected from the previous results for polygons. Then we further train this model on 4800 3d data points. It turns out that the model originally trained for 2d data has successfully learn the 3d data. In fact, this result is just as good as training 3d data only with fixed length 10 as in Table \ref{plk10_vol3d}. See Table \ref{vol2dto3d}.
\begin{table}[h]
\centering
\begin{tabular}{|c|c|c|c|c|c|}
\hline
MAE & \multicolumn{5}{c|}{0.301} \\ \hline\hline
Error & $\pm0.5$ & $\pm1$ & $\pm2$ & $\pm3$ & $\pm4$ \\
Accuracy & 0.825 & 0.977 & 1.000 & 1.000 & 1.000 \\
\hline
\end{tabular}
\caption{The MAE and accuracies for predicting polytope volumes with input of fixed length 10 using a model pre-trained for 2d data.}\label{vol2dto3d}
\end{table}

As a sanity check, we now train the model only with 2d data and extrapolate it to 3d data without further fitting any 3d samples. We find that even at $\text{err}=\pm5$, the accuracy can only reach 0.007 as one may expect.

Since the results are always good for Pl\"ucker coordinates of length 10 no matter whether we have pre-trained 2d data or not, one may wonder whether pre-training 2d data would make a difference for varying input lengths. Since the longest input vector would have length 560 (from 3d data), we pad all the 2d and 3d data with 0's to length 560. As padding many 0's would introduce noise, we now include more data for training for both 2d and 3d polytopes. Same as the steps above, we first train the model with 76000 2d samples and then further fit the model with 45000 3d samples. The results are reported in Table \ref{vol2dto3dvarying}.
\begin{table}[h]
\centering
\begin{tabular}{|c|c|c|c|c|c|}
\hline
MAE & \multicolumn{5}{c|}{2.149} \\ \hline\hline
Error & $\pm0.5$ & $\pm1$ & $\pm2$ & $\pm3$ & $\pm4$ \\
Accuracy & 0.152 & 0.302 & 0.570 & 0.751 & 0.860 \\
\hline
\end{tabular}
\caption{The MAE and accuracies for predicting the polytope volumes with input of varying lengths (30000 samples) using the model pre-trained \textbf{with} 2d data (15000 samples).}\label{vol2dto3dvarying}
\end{table}
\begin{table}[h]
\centering
\begin{tabular}{|c|c|c|c|c|c|}
\hline
MAE & \multicolumn{5}{c|}{2.193} \\ \hline\hline
Error & $\pm0.5$ & $\pm1$ & $\pm2$ & $\pm3$ & $\pm4$ \\
Accuracy & 0.148 & 0.295 & 0.550 & 0.748 & 0.861 \\
\hline
\end{tabular}
\caption{The MAE and accuracies for predicting the polytope volumes with input of varying lengths (30000 samples) \textbf{without} pre-training any 2d data.}\label{vol3dvaryingcompare}
\end{table}

Again, the model can still learn the 3d data pretty well. However, if we only train same amount of 3d data (45000 samples) without pre-training 2d data, the accuracy would not make a big difference as shown in Table \ref{vol3dvaryingcompare}.

We can also perform the similar test for dual volumes. Again, we find that a model pre-trained with 2d data gives similar performance as a model without pre-training does. Although the polytope properties for 2d polygons and 3d polytopes in terms of Pl\"ucker coordinates might share some common ground, the above results show that different dimensions would still have certain different features.

\subsection{The MDS Embedding}\label{mds}
We have used different models to predict various properties of the polytopes with Pl\"ucker coordinates. In particular, the performance of predicting volumes is most impressive amongst them. In fact, the volume $V$ is the sum of a subset of Pl\"ucker coordinates in the input vector, that is,
\begin{equation}
    V=\sum_{\alpha\in\mathcal{V}}p^\alpha\label{subsetsum}
\end{equation}
where $(p^0,p^1,\dots,p^{l-1})$ is the Pl\"ucker coordinates (of length $l$) and $\mathcal{V}\subset\{0,1,\dots,l-1\}$.

Now let us apply manifold embedding with the help of \texttt{scikit-learn} \cite{scikit-learn} and \texttt{Yellowbrick} \cite{bengfort_yellowbrick_2019} to understand how the machine learning models make the predictions and what results we can extract from such analysis.

\begin{figure}[h]
	\centering
	\begin{subfigure}{6cm}
	\includegraphics[width=6cm]{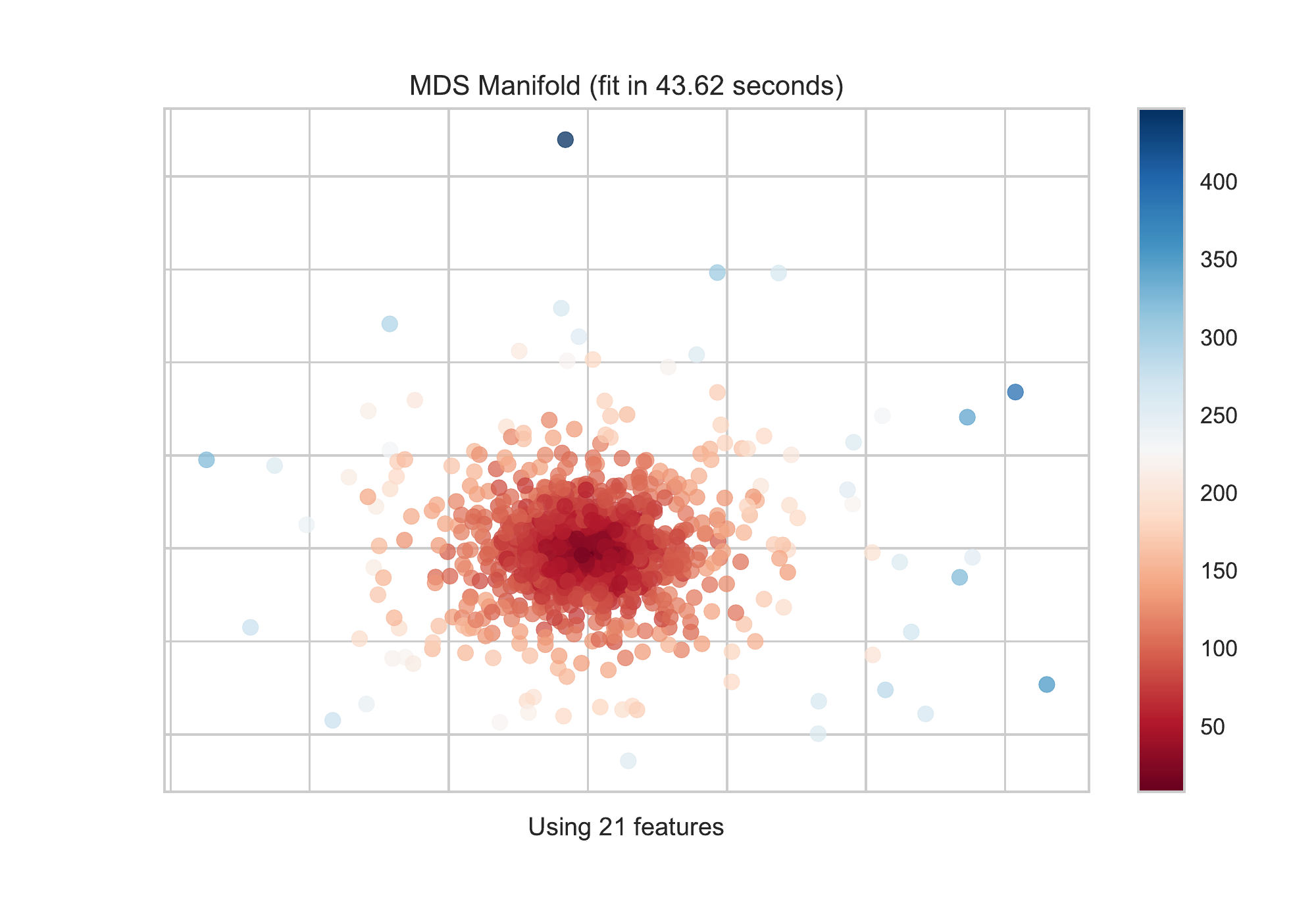}
    \caption{}
	\end{subfigure}
	\begin{subfigure}{6cm}
	\includegraphics[width=6cm]{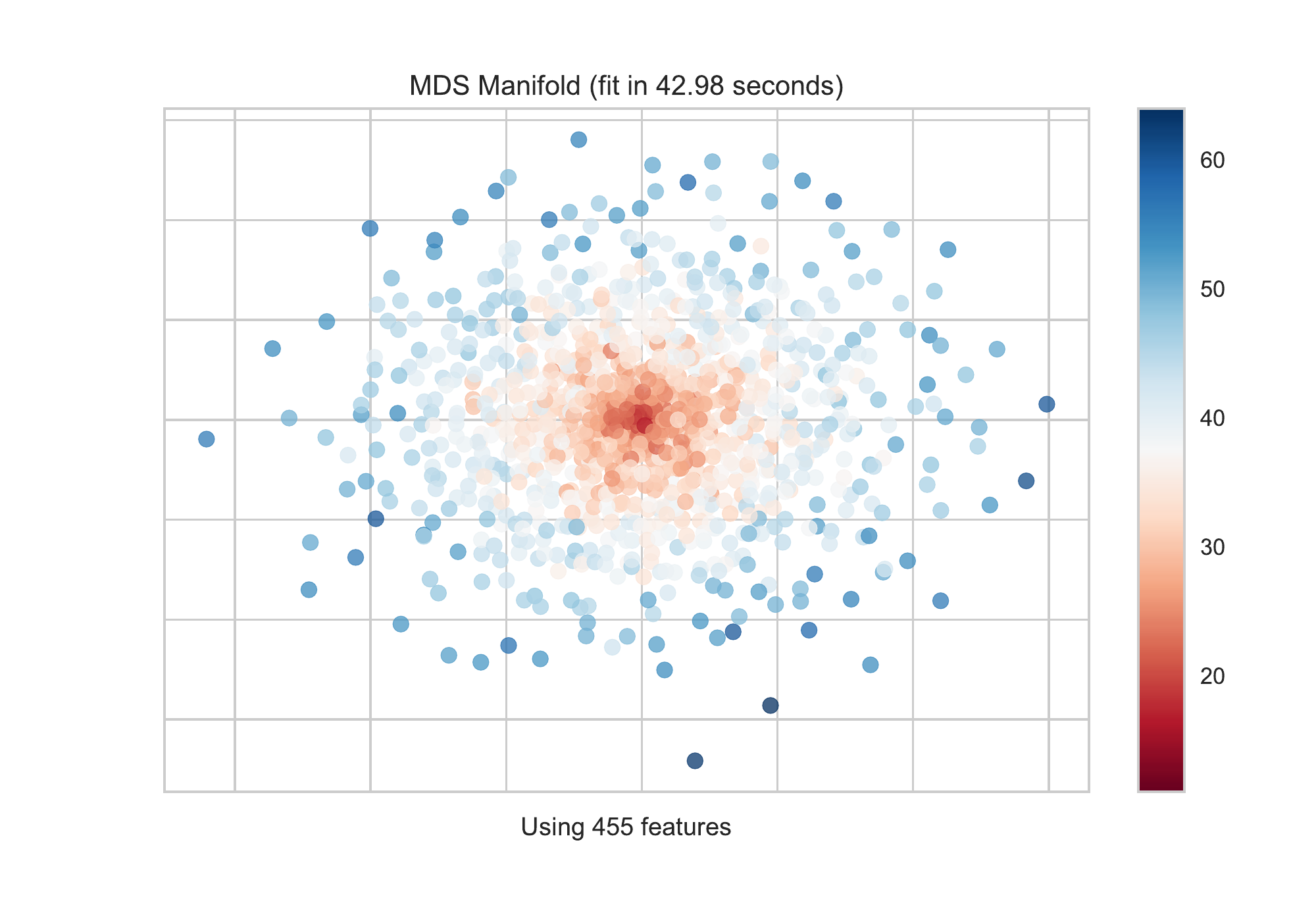}
    \caption{}
	\end{subfigure}
	\caption{The (2-component) MDS projections of Pl\"ucker coordinates for (a) polygons and (b) polytopes. Different colours of data points indicate different \textit{volumes}. Here, we use 1000 random samples for each plot as an illustration.}\label{volmds}
\end{figure}

It turns out that multi-dimensional scaling (MDS) \cite{borg2005modern} leads to nice distributions of the data\footnote{We have also tried PCA and other different manifold embeddings. However, only MDS could give scattering plots with such distributions.}. See Figure \ref{volmds}.

In general, the MDS applied here projects a Pl\"ucker vector $\bm{p}=(p^0,\dots,p^{l-1})$ to a $k$-dimensional vector $\bm{x}=(x^0,\dots,x^{k-1})$ ($k<l$) by minimizing the cost function known as the stress:
\begin{equation}
    S=\left(\sum_{i<j\leq N}(D_{ij}-d_{ij})^2\right)^{1/2},
\end{equation}
where $N$ is the number of samples and $D_{ij}$ ($d_{ij}$) denotes the Euclidean distance between $\bm{p}_i$ ($\bm{x}_i$) and $\bm{p}_j$ ($\bm{x}_j$).

\begin{figure}[h]
	\centering
	\begin{subfigure}{6cm}
	\includegraphics[width=6cm]{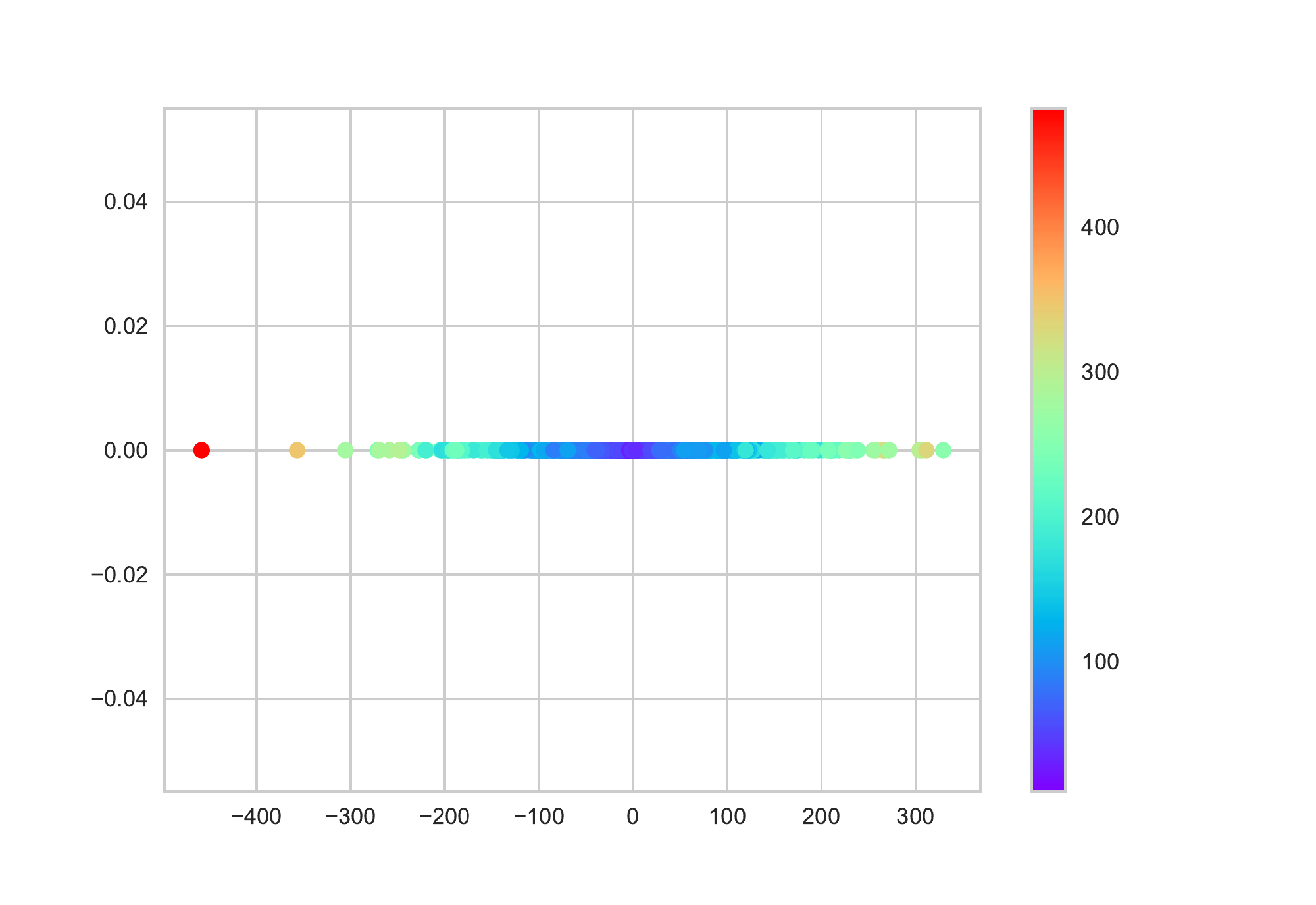}
    \caption{}
	\end{subfigure}
	\begin{subfigure}{6cm}
	\includegraphics[width=6cm]{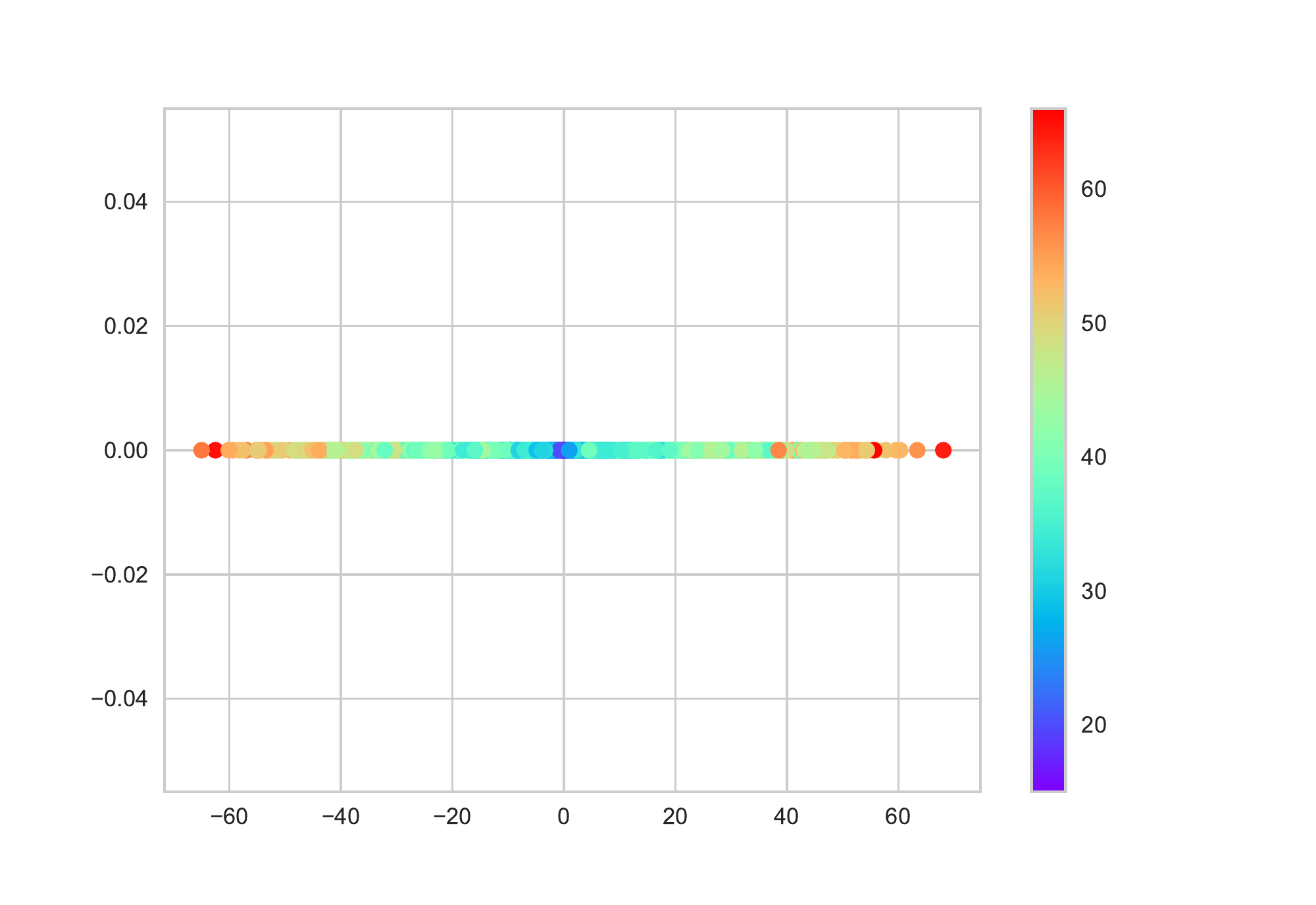}
    \caption{}
	\end{subfigure}
	\caption{The \emph{$1$-component} MDS projections of Pl\"ucker coordinates for (a) polygons and (b) polytopes. Different colours of data points indicate different volumes. Here, we use 1000 random samples for each plot as an illustration.}\label{volmds1component}
\end{figure}

From Figure \ref{volmds}, we can see that the vectors $\bm{x}$ are distributed around a centre point. The closer a vector is to the centre, the smaller volume the corresponding polytope (either 2d or 3d) would have. Since this is a mathematical fact, the number of samples should not affect the results and hence minimizing $S$ is equivalent to the goal of finding $\bm{x}$ such that $D_{ij}-d_{ij}\approx0$ (for any $i,j$). Since the positions of the data points with different colours/volumes mainly depend on the distances/radii from the origin in Figure \ref{volmds}, we may also plot a 1-component MDS projection, i.e., $\bm{x}=(x^0)\equiv x$. As shown in Figure \ref{volmds1component}, the points with different colours/volumes are still determined by their distances to the origin.

This seems to hint some possible relation between Pl\"ucker coordinates and polytope volumes.
\begin{conjecture}
There exists a map $f:\bigcup\limits_nG(n-d,n)\rightarrow\mathbb{R},\,\bm{p}\mapsto x$ such that
\begin{equation}
    V\sim|x|+1,
\end{equation}
where $V$ is the corresponding polytope volume of the Pl\"ucker coordinates $\bm{p}$, and $1$ indicates the smallest (normalized) volume among the polytopes.\label{conjvol}
\end{conjecture}
Recall that $G(n-d,n)$ is the Grassmannian which parametrizes the Pl\"ucker coordinates of $n$-topes in $d$ dimensions. Since the Pl\"ucker vectors would have different lengths, we pad 0's to make $\bm{p}$ have the same length $l$.

In terms of minimizing $S$, for any two Pl\"ucker vectors $\bm{p}$ and $\bm{q}$ with $x=f(\bm{p})$ and $y=f(\bm{q})$,
\begin{equation}
    D(\bm{p},\bm{q})-d(x,y)=0.
\end{equation}
Hence,
\begin{equation}
    (|x|-|y|)^2=x^2+y^2-2|xy|=\sum_{\alpha=1}^l((p^{\alpha})^2+(q^{\alpha})^2-2p^{\alpha}q^{\alpha}).
\end{equation}
Since $x=f(\bm{p})$ and $y=f(\bm{q})$, the crossing terms $p^{\alpha}q^{\alpha}$ on the right hand side should only come from the crossing term $|xy|$ on the left hand side. Therefore, by \eqref{subsetsum}, we have
\begin{equation}
\begin{split}
    &x^2=\left(\sum_{\alpha\in\mathcal{V}_{\bm{p}}}p^{\alpha}\right)^2=\sum_{\alpha\in\mathcal{V}_{\bm{p}}}(p^{\alpha})^2+\sum_{\alpha\neq\beta\in\mathcal{V}_{\bm{p}}}p^{\alpha}p^{\beta},\\
    &|xy|=\sum_{\alpha=1}^lp^{\alpha}q^{\alpha}+\sum_{\alpha\neq\beta\in\mathcal{V}_{\bm{p}}}p^{\alpha}p^{\beta}+\sum_{\alpha\neq\beta\in\mathcal{V}_{\bm{q}}}q^{\alpha}q^{\beta}-\sum_{\alpha\neq\beta\not\in\mathcal{V}_{\bm{p}}}p^{\alpha}p^{\beta}-\sum_{\alpha\neq\beta\not\in\mathcal{V}_{\bm{q}}}q^{\alpha}q^{\beta}.
\end{split}
\end{equation}

Overall, the MDS embedding hints some non-trivial relation between Pl\"ucker coordinates and volumes as in Figure \ref{volmds} and Figure \ref{volmds1component}. It would also be interesting to study its connection to the Pl\"ucker relations.

\begin{figure}[h]
	\centering
	\begin{subfigure}{6cm}
	\includegraphics[width=6cm]{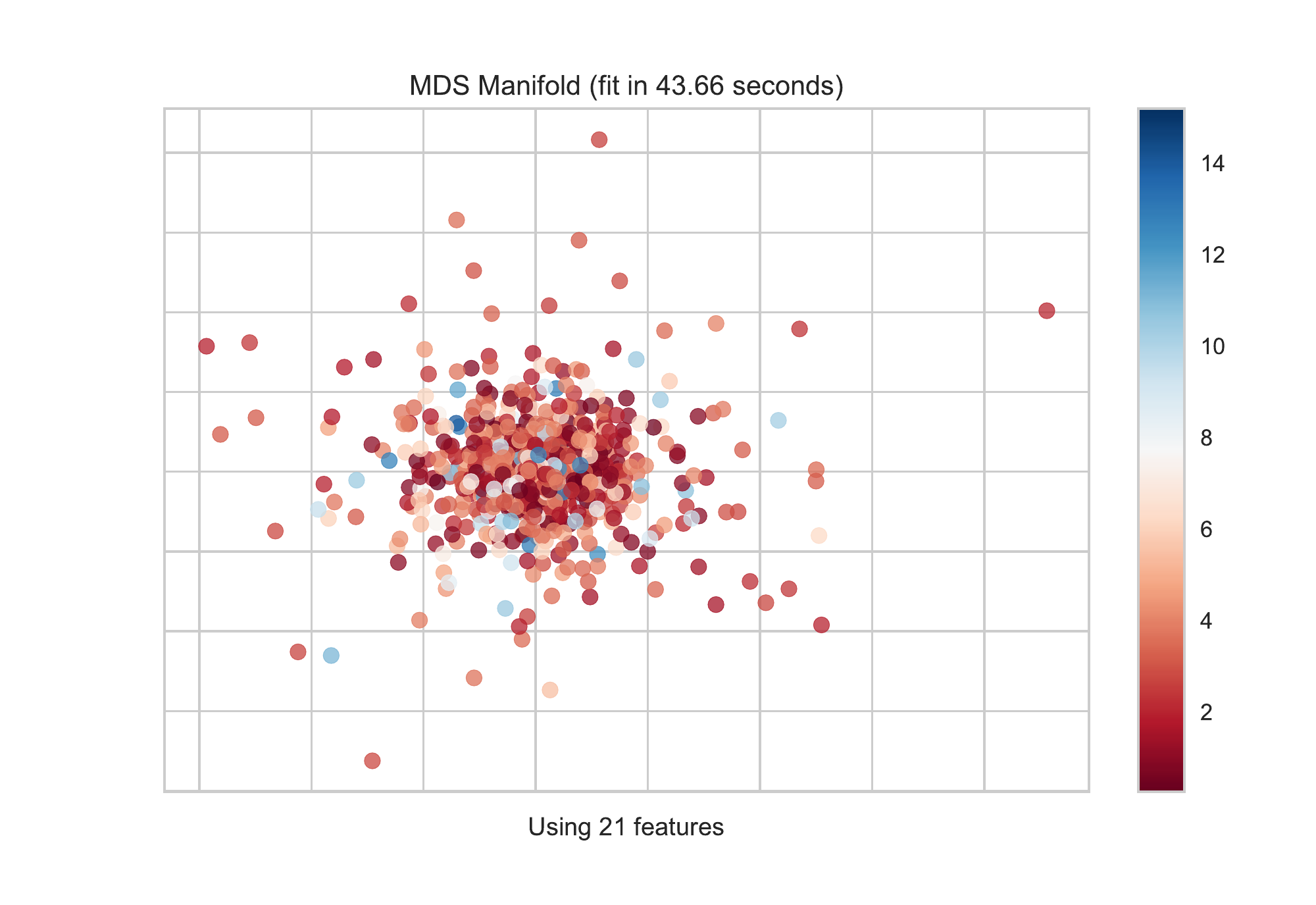}
    \caption{}
	\end{subfigure}
	\begin{subfigure}{6cm}
	\includegraphics[width=6cm]{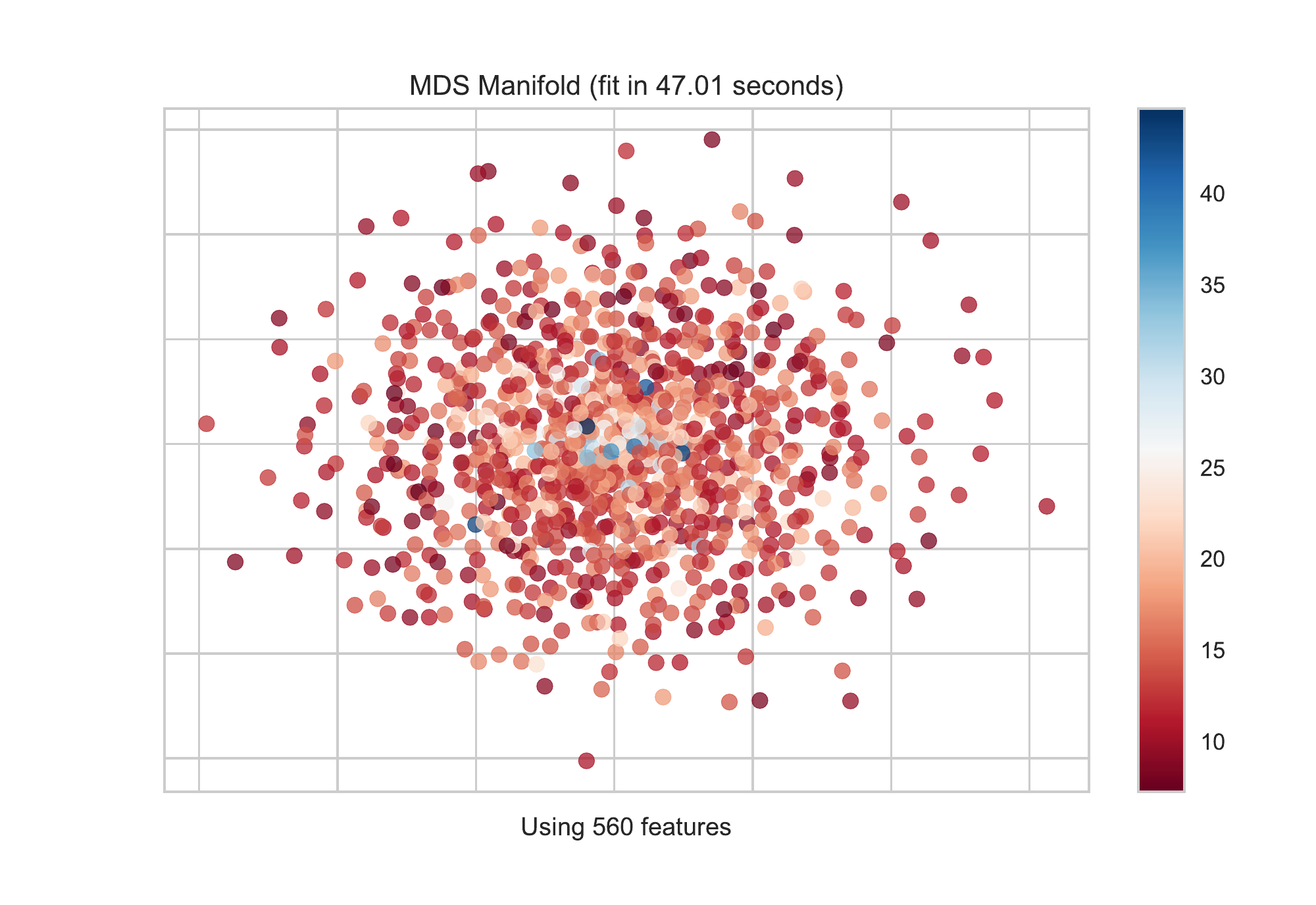}
    \caption{}
	\end{subfigure}
	\caption{The MDS projections of Pl\"ucker coordinates for (a) polygons and (b) polytopes. Different colours of data points indicate different \textit{dual volumes}. Here, we use 1000 random samples for each plot as an illustration.}\label{dualvolmds}
\end{figure}

\paragraph{A comment on dual volumes} We may also apply MDS to dual volumes. As shown in Figure \ref{dualvolmds}, although data points with larger dual volumes tend to lie closer to the centres in the plots (as opposed to volumes), the ascending of dual volumes from centres to peripheries is not comparable to the case for volumes.

\section{Summary \& Outlook}
In this paper, we explored the machine learning for polytopes in terms of their Pl\"ucker coordinates. Here, we summarize the main results in Table \ref{summarytable}.

These results indicate that it could be possible to find certain relations, albeit very non-trivial, between Pl\"ucker coordinates and various properties of the polytopes. In particular, some of the results show much better performance using Pl\"ucker coordinates than using the vertices. From Figure \ref{volmds1component}, we also found that the volumes have nice distributions under MDS embedding for Pl\"ucker coordinates. This indicates possible relations between volumes and Pl\"ucker coordinates. It is then natural to ask what the exact relation is and whether this could be connected to our current knowledge on Pl\"ucker coordinates such as Pl\"ucker relations etc. Moreover, it would be important to see how neural network would perform when predicting the MDS-projected data from Pl\"ucker coordinates as well as predicting volumes from such data \cite{MLKS}.

It is also natural to expect that the above study could be extended to polytopes in other dimensions. For instance, in 4d, the Kreuzer-Skarke (KS) programme of studying Calabi-Yau manifolds as hypersurfaces (and, more generally, complete intersections) in toric varieties from reflexive polytopes is still a very active area of research \cite{Braun:2012vh,Altman:2014bfa,He:2015fif,He:2017gam,Altman:2018zlc,Demirtas:2018akl,Demirtas:2020dbm}. The important role played by the KS database in both mathematics and physics makes it a particularly interesting future direction for this research.

\begin{table}[h]
\centering
\begin{tabular}{V{4}cV{4}c|c|c|cV{4}}
\hlineB{4}
\begin{tabular}[c]{@{}c@{}}Polygons\\ (frequency)\end{tabular}   & \begin{tabular}[c]{@{}c@{}}Triangles\\ (277)\end{tabular} & \begin{tabular}[c]{@{}c@{}}Quadrilaterals\\ (7041)\end{tabular} & \begin{tabular}[c]{@{}c@{}}Pentagons\\ (16637)\end{tabular} & \begin{tabular}[c]{@{}c@{}}Hexagons\\ (3003)\end{tabular} \\ \hlineB{4}
Input & \multicolumn{4}{cV{4}}{Vertices} \\ \hline
Output & \multicolumn{4}{cV{4}}{Volume} \\ \hline
MAE & 4.941 & 10.012 & 8.640 & 14.947 \\ \hline
Accuracy & 0.945 & 0.945 & 0.926 & 0.826 \\ \hlineB{4}
Input & \multicolumn{4}{cV{4}}{Pl\"ucker coordinates} \\ \hline
Output & \multicolumn{4}{cV{4}}{Volume} \\ \hline
MAE & 0.209 & 0.625 & 1.051 & 3.359 \\ \hline
Accuracy & 1.000 & 1.000 & 1.000 & 0.997  \\ \hlineB{4}
Input & \multicolumn{4}{cV{4}}{Pl\"ucker coordinates} \\ \hline
Output  & \multicolumn{4}{cV{4}}{Dual Volume} \\ \hline
MAE & 1.181 & 0.642 & 0.818 & 0.941 \\ \hline
Accuracy & 0.501 & 0.754 & 0.638 & 0.557 \\ \hlineB{4}
Input & \multicolumn{4}{cV{4}}{($\ell$-1)-Pl\"ucker coordinates \& Volume} \\ \hline
Output & \multicolumn{4}{cV{4}}{The missing Pl\"ucker coordinate} \\ \hline
MAE & 0.295 & 2.599 & 2.061 & 3.169 \\ \hline
Accuracy & 1.000 & 0.983 & 0.994 & 0.971 \\ \hlineB{4}
\begin{tabular}[c]{@{}c@{}}Polytopes\\ (frequency)\end{tabular}  & \multicolumn{4}{cV{4}}{\begin{tabular}[c]{@{}c@{}}Toric canonical Fano polytopes\\ (780000)\end{tabular}} \\ \hline
Input & \multicolumn{4}{cV{4}}{Pl\"ucker coordinates} \\ \hline
Output & Volume & Dual Volume & \multicolumn{2}{cV{4}}{Reflexivity} \\ \hline
MAE & 1.680 & 2.590 & \multicolumn{2}{cV{4}}{N/A} \\ \hline
Accuracy & 0.936 & 0.890 & \multicolumn{2}{cV{4}}{0.813} \\ \hlineB{4}
\end{tabular}
\caption{Summary of the key machine learning results. Notice that for polygons with Pl\"ucker coordinates, the ``Accuracy'' in the table is the accuracy ($\pm0.05\times\text{range}$). Likewise, for those involving polytopes with Pl\"ucker coordinates, the ``Accuracy'' in the table is the accuracy ($\pm4$). For reflexivity the result is simply the accuracy itself since this is a binary classification. Moreover, the input used for classifying reflexivity is one-hot encoding of Pl\"ucker coordinates for the result listed here.}
\label{summarytable}
\end{table}

\section*{Acknowledgement}
JB~is supported by a~CSC scholarship.
YHH~would like to thank~STFC for grant ST/J00037X/1.
EH~would like to thank~STFC for a~PhD studentship.
JH~is supported by a Nottingham Research Fellowship.
AK~is supported by~EPSRC Fellowship~EP/N022513/1.
SM~is funded by a~SMCSE Doctoral Studentship. This collaboration was made possible by a Focused Research Workshop grant from the Heilbronn Institute for Mathematical Research.

\addcontentsline{toc}{section}{References}
\bibliographystyle{utphys}
\bibliography{references}

\end{document}